# IMPLICIT-EXPLICIT FORMULATIONS FOR DISCONTINUOUS GALERKIN NON-HYDROSTATIC ATMOSPHERIC MODELS

**Sohail Reddy** *
Department of Applied Mathematics
Naval Postgraduate School
Monterey, CA
sredd001@fiu.edu

**Maciej Waruszewski**
Sandia National Laboratories
Albuquerque, NM
mwarusz@sandia.gov

**Felipe A. V. de Braganca Alves**
Department of Applied Mathematics
Naval Postgraduate School
Monterey, CA
felipe.alves.br@nps.edu

**Francis X. Giraldo**
Department of Applied Mathematics
Naval Postgraduate School
Monterey, CA
fxgirald@nps.edu

## ABSTRACT

This work presents IMplicit-EXplicit (IMEX) formulations for discontinuous Galerkin (DG) discretizations of the compressible Euler equations governing non-hydrostatic atmospheric flows. In particular, we show two different IMEX formulations that not only treat the stiffness due to the governing dynamics but also the domain discretization. We present these formulations for two different equation sets typically employed in atmospheric modeling. For both equation sets, efficient Schur complements are derived and the challenges and remedies for deriving them are discussed. The performance of these IMEX formulations of different orders are investigated on both 2D (box) and 3D (sphere) test problems and shown to achieve their theoretical rates of convergence and their efficiency with respect to both mesoscale and global applications are presented.

## 1 Introduction

The discontinuous Galerkin (DG) method has been used successfully for developing compressible Euler and Navier-Stokes models. However, one of the main drawbacks of the method is that it is not competitive with either low-order or continuous (CG) methods due to the smaller time-step required for stability. It is well-known that high-order methods impose a more stringent Courant-Friedrichs-Lewy (CFL) condition compared to low-order methods (see, e.g., [1]). It is less clearly understood that most DG methods in the literature impose a more stringent CFL condition than CG methods; the reason is twofold: (1) upwind-based numerical fluxes move the eigenvalues from the imaginary to the real axis and so the time-integration must be stable for such eigenspectra (e.g., Runge-Kutta methods satisfy this condition but neutrally stable methods do not); (2) most Godunov methods in use today (e.g., FV and DG methods) rely on one-dimensional numerical fluxes - this means that if the flow moves diagonally across a gridpoint then it must cross two element edges. If the numerical method is not constructed in such a manner (as in the evolutionary DG method described in [2]) then smaller time-steps are required. However, our experience has shown that using standard (one-dimensional fluxes that require a smaller time-step) give faster time-to-solution than more complicated fluxes such as those used in [2].

*Also with Center for Applied Scientific Computing (CASC) at Lawrence Livermore National Laboratory, Livermore, CA, USA. (reddy6@llnl.gov)

In order to circumvent the stringent CFL condition imposed by explicit methods, researchers have explored implicit-explicit (IMEX) methods. The first two such works found in the literature are by [3] and [4] where the authors begin with standard double Butcher tableaux whereby the stiff portion is handled implicitly while the non-stiff part is handled explicitly. In both works they begin with nonlinear handling of the stiff implicit part but in [3], the authors linearize the problem relying on a specific numerical flux while those in [4] maintain the nonlinear structure and use Newton's method to linearize the problem. The work that we describe here continues the work first proposed in [5, 6] where we rely on the linearization of the implicit nonlinear problem using a reference solution, typical of semi-implicit methods [7, 8] used in the atmospheric modeling community. This approach has been recently analyzed in the so-called RS-IMEX method (the RS denotes the *reference solution* linearization) in [9–13]. Our goal is to complete the work that we began in [6] where we sought a Schur complement to the construction of the DG IMEX method but were only able to do so for a specific class of boundary conditions (namely, periodic in all directions). In what follows, we generalize this approach to reflection boundary conditions, although other types of boundary conditions are also possible. Furthermore, we show that the eigenvalues of the resulting pseudo-Helmholtz operator resulting from the Schur complement are all real and positive, indicating that the resulting operator is symmetric positive-definite, similar to its CG counterpart (see, e.g., [5, 14]). We then construct Schur complements for the full 3D IMEX problem (when the grid is isotropic) as well as for the 1D IMEX problem (when the stiffness results from one direction, as occurs along the vertical direction in geophysical fluid dynamics applications of atmosphere/climate models).

## 2 Governing Equations

This work presents semi-implicit formulations for two different equation sets often used in non-hydrostatic atmospheric modeling. We limit our discussions and formulations to equations sets in conservation form. The nomenclature that we use below was first defined in [15] where various forms of the Euler equations were defined, discretized, and compared. We only choose to discuss the equations in conservation form since these are the only practical options for using the discontinuous Galerkin method. We consider balance laws written in the form

$$\frac{\partial \mathbf{q}}{\partial t} + \nabla \cdot \mathbf{F} = \mathbf{S} \tag{2.1}$$

where $\mathbf{q}$, $\mathbf{F}$, and $\mathbf{S}$ are vectors of prognostic, conserved variables and their total fluxes and sources, respectively. These vectors are defined in the subsequent sections.

### 2.1 Equation Set2C

This equation set is used in its conservation form in [16, 17] and in non-conservation form in [14]. This equation set features the standard formulations for the continuity and momentum equations but utilizes the thermodynamic equation in terms of potential temperature. This can be advantageous in atmospheric modeling as this quantity is unaffected by the rising and sinking motion of the fluid over obstacles (such as mountains) or in large-scale atmospheric turbulence; the potential temperature is also a useful measure of the static stability of the atmosphere [18, Ch. 8]. One drawback of this equation set is that its derivation (from the first law of thermodynamics) assumes a constant composition atmosphere (i.e., gas constants remain constant); therefore, this equation set is not a good choice for atmospheric simulations reaching into the thermosphere. Another drawback (shown in [15]) is that the equation for computing pressure is expensive due to the required exponentials. This equation set is written as Eq. (2.1) with

$$\mathbf{q}_{2C} = \begin{Bmatrix} \rho \\ \mathbf{U} \\ \Theta \end{Bmatrix}, \quad \mathbf{F}_{2C} = \begin{Bmatrix} \mathbf{U} \\ \dfrac{\mathbf{U} \otimes \mathbf{U}}{\rho} + P\mathbf{I} \\ \left(\dfrac{\Theta}{\rho}\right) \mathbf{U} \end{Bmatrix}, \quad \mathbf{S}_{2C} = \begin{Bmatrix} 0 \\ -\left(\rho \nabla \phi + 2\Omega \times \mathbf{U}\right) \\ 0 \end{Bmatrix} \tag{2.2}$$

where the conserved, prognostic variables are $(\rho, \mathbf{U}, \Theta)$, $\mathbf{U} = \rho\mathbf{u} = (\rho u, \rho v, \rho w)^{\mathcal{T}}$, $\Theta = \rho\theta$ is the potential temperature density, $\phi = gr$ is the geopotential height with $g$ representing the gravitational constant and $r$ representing the height coordinate (in flow-in-a-box, it is the Cartesian z-component and on the sphere it represents the radial coordinate), $\mathbf{I}$ is a rank-3 identity tensor, and $^{\mathcal{T}}$ is the transpose operator. The pressure $P$ is given as

$$P = P_A \left(\frac{R\Theta}{P_A}\right)^{\gamma} \tag{2.3}$$

where $P_A$ is the reference pressure at the surface ($P_A = 1\times10^5$ Pa), $R = c_p - c_v$ is the specific gas constant, $c_p$ and $c_v$ are the specific heats at constant pressure and volume, respectively, and $\gamma = c_p/c_v$ is the ratio of specific heats. Introducing splitting about some reference field for density, $\rho(\mathbf{x},t) = \rho_0(r) + \rho'(\mathbf{x},t)$, pressure, $P(\mathbf{x},t) = P_0(r) + P'(\mathbf{x},t)$, and potential temperature, $\Theta(\mathbf{x},t) = \Theta_0(r) + \Theta'(\mathbf{x},t)$ where the reference values $(\rho_0, P_0, \Theta_0)$ are in hydrostatic balance ($\nabla P_0 = -\rho_0 \nabla\phi$), Eq. (2.2) can be written in pertubation form as Eq. (2.1) with $\mathbf{q} := \mathbf{q}'_{2C}$, $\mathbf{F} := \mathbf{F}'_{2C}$ and, $\mathbf{S} := \mathbf{S}'_{2C}$ defined as

$$\mathbf{q}'_{2C} = \begin{Bmatrix} \rho' \\ \mathbf{U} \\ \Theta' \end{Bmatrix}, \quad \mathbf{F}'_{2C} = \begin{Bmatrix} \mathbf{U} \\ \dfrac{\mathbf{U}\otimes\mathbf{U}}{\rho} + P'\mathbf{I} \\ \left(\dfrac{\Theta}{\rho}\right)\mathbf{U} \end{Bmatrix}, \quad \mathbf{S}'_{2C} = \begin{Bmatrix} 0 \\ -\left(\rho'\nabla\phi + 2\Omega\times\mathbf{U}\right) \\ 0 \end{Bmatrix}. \tag{2.4}$$

## 2.2 Equation Set3C

This equation set is the most common form of the Euler equations in CFD [19] and has been used in two recent dynamical cores related to the CliMa model [20,21]. This equation set shares the continuity and momentum equations with Set2C but differs in the thermodynamic equation and the equation of state. One of the advantages of this equation set is that it allows for conservation of not only mass and momentum, but also total energy. Moreover, since it is the equation of choice in CFD, it has been used in many recent works that rigorously prove the stability of this equation set using DG methods even at high-order (see [21,22]). Yet another advantage of this equation set is that its derivation does not assume anything about the composition of the atmosphere and so is quite general. A disadvantage of this equation set, at least in the context of atmospheric modeling, is that much of the physical parameterizations require either potential temperature or temperature as inputs and so one would need to extract this variable from total energy. One other likely disadvantage (based on our experience) is that this equation set requires more care to stabilize since it is completely conservative (i.e., non-dissipative) and so the numerical method has to be designed carefully to maintain stability. Set3C is written as Eq. (2.1) with

$$\mathbf{q}_{3C} = \begin{Bmatrix} \rho \\ \mathbf{U} \\ E \end{Bmatrix}, \quad \mathbf{F}_{3C} = \begin{Bmatrix} \mathbf{U} \\ \dfrac{\mathbf{U}\otimes\mathbf{U}}{\rho} + P\mathbf{I} \\ \dfrac{(E+P)}{\rho}\mathbf{U} \end{Bmatrix}, \quad \mathbf{S}_{3C} = \begin{Bmatrix} 0 \\ -\left(\rho\nabla\phi + 2\Omega\times\mathbf{U}\right) \\ 0 \end{Bmatrix} \tag{2.5}$$

where the conserved, prognostic variables are $(\rho, \mathbf{U}, E)$, $E$ is the total energy given by $E = \rho c_v T + \dfrac{\mathbf{U}\cdot\mathbf{U}}{2\rho} + \rho\phi$ which denote the internal, kinetic, and potential energies, respectively. The pressure $P$ is given as

$$P = \rho R T \equiv (\gamma - 1)\left(E - \frac{\mathbf{U}\cdot\mathbf{U}}{2\rho} - \rho\phi\right). \tag{2.6}$$

Note that Eq. (2.6) contains no exponentials, as does Set2C, and is therefore cheaper to compute. Since this equation uses the total energy as the thermodynamic variable and is written in conservation form, we can conserve energy; this is not guaranteed for Set2C.

Splitting the density and pressure as done for Set2C and total energy as $E(\mathbf{x},t) = E_0(r) + E'(\mathbf{x},t)$, and assuming the reference fields are in hydrostatic balance, Set3C can be written in pertubation form as Eq. (2.1) with $\mathbf{q} := \mathbf{q}'_{3C}$, $\mathbf{F} := \mathbf{F}'_{3C}$ and, $\mathbf{S} := \mathbf{S}'_{3C}$ defined as

$$\mathbf{q}'_{3C} = \begin{Bmatrix} \rho' \\ \mathbf{U} \\ E' \end{Bmatrix}, \quad \mathbf{F}'_{3C} = \begin{Bmatrix} \mathbf{U} \\ \dfrac{\mathbf{U}\otimes\mathbf{U}}{\rho} + P'\mathbf{I} \\ \dfrac{(E+P)}{\rho}\mathbf{U} \end{Bmatrix}, \quad \mathbf{S}'_{3C} = \begin{Bmatrix} 0 \\ -\left(\rho'\nabla\phi + 2\Omega\times\mathbf{U}\right) \\ 0 \end{Bmatrix}. \tag{2.7}$$

# 3 IMplicit-EXplicit Formulations

## 3.1 Semi-Implicit Time Integration

The governing equations for both Set2C and Set3C can be written in vector form as

$$\begin{aligned} \frac{\partial \mathbf{q}}{\partial t} &= S(\mathbf{q}) = \{E(\mathbf{q})\} + [I(\mathbf{q})] \\ &= \{S(\mathbf{q}) - L(\mathbf{q})\} + [L(\mathbf{q})] \end{aligned} \tag{3.1}$$

where for Set2C and Set3C, $\mathbf{q}_{2C} := (\rho, \mathbf{U}, \Theta)$ and $\mathbf{q}_{3C} := (\rho, \mathbf{U}, E)$, respectively, $S(\mathbf{q})$ contains all spatial terms, $E(\mathbf{q})$ contains the operators to be treated explicitly, $I(\mathbf{q})$ contains the operators to be treated implicitly, and $L(\mathbf{q})$ is a linear operator to be defined in the following sections. The semi-implicit formulation allows the full problem $(S(\mathbf{q}))$ to be decomposed into subsets that can be treated differently. In Eq. (3.1), the terms in the curly brackets $\{\cdot\}$ contain the slow moving advection terms and are treated explicitly, whereas the terms in the square brackets $[\cdot]$ contain the fast moving acoustic and gravity terms and are treated implicitly. This implicit treatment of the restrictive dynamics removes the CFL constraint due to the fast propagating acoustic and gravity waves and allows for larger time steps that are dictated only by the slow dynamics. The implicit treatment of linear operators requires the solution of a linear system whereas the treatment of nonlinear operators requires a nonlinear system to be solved.

Consider a single-step multistage method written as

$$\mathbf{Q}^{(i)} = \mathbf{q}^n + \Delta t \left[\sum_{j=1}^{i-1} a_{ij} \left(S\left(\mathbf{Q}^{(j)}\right) - L\left(\mathbf{Q}^{(j)}\right)\right)\right] + \Delta t \left[\sum_{j=1}^{i-1} \widetilde{a}_{ij} \left(L\left(\mathbf{Q}^{(j)}\right)\right)\right] + \Delta t \widetilde{a}_{ii} L\left(\mathbf{Q}^{(i)}\right) \tag{3.2}$$

where $\mathbf{Q}^{(i)}$ is the $i^{th}$ stage approximation, $\Delta t$ is the time-step, $a_{ij}$, $\widetilde{a}_{ij}$, $b_i$, and $\widetilde{b}_i$ are the coefficients in the Butcher tableaux for the explicit and implicit components, respectively. Following [14] we define

$$\mathbf{q}_{tt} = \mathbf{Q}^{(i)} + \sum_{j=1}^{i-1} \frac{\widetilde{a}_{ij} - a_{ij}}{\widetilde{a}_{ii}} \mathbf{Q}^{(j)} \tag{3.3}$$

and $\widehat{\mathbf{q}}$ as

$$\widehat{\mathbf{q}} = \mathbf{q}^n + \sum_{j=1}^{i-1} \frac{\widetilde{a}_{ij} - a_{ij}}{\widetilde{a}_{ii}} \mathbf{Q}^{(j)} + \Delta t \sum_{j=1}^{i-1} a_{ij} S\left(\mathbf{Q}^{(j)}\right). \tag{3.4}$$

Eq. (3.2) can be written compactly as

$$\mathbf{q}_{tt} = \widehat{\mathbf{q}} + \alpha L(\mathbf{q}_{tt}) \tag{3.5}$$

where $\alpha = \Delta t \widetilde{a}_{ii}$. The stage value $\mathbf{Q}^{(i)}$ can be obtained from $\mathbf{q}_{tt}$ using Eq. (3.3) and the solution update at $t_{n+1}$ can be constructed as

$$\mathbf{q}^{n+1} = \mathbf{q}^n + \Delta t \sum_{i=1}^{s} b_i \left[S\left(\mathbf{Q}^{(i)}\right) - L\left(\mathbf{Q}^{(i)}\right)\right] + \Delta t \sum_{i=1}^{s} \widetilde{b}_i L\left(\mathbf{Q}^{(i)}\right) \tag{3.6}$$

Assuming that the coefficients in the Butcher tableaux satisfy $b_i = \widetilde{b}_i$ needed for conservation (see [14]), we obtain

$$\mathbf{q}^{n+1} = \mathbf{q}^n + \Delta t \sum_{i=1}^{s} b_i S\left(\mathbf{Q}^{(i)}\right). \tag{3.7}$$

For the IMEX scheme to be consistent on the continuous and discrete levels, the implicitly and explicitly treated components should be separable on both levels. We discuss the necessary conditions for this in Appendix C. In what follows, we derive the *No-Schur* and Schur forms for equation sets 2C and 3C, however, to simplify the exposition, we perform the derivation using continuous differential operators. In the appendices, we substitute the continuous differential operators with their discrete forms to complete the discussion.

### 3.2 IMplicit-EXplicit Formulation: Set2C

#### 3.2.1 3D-IMEX Formulation

Let us now describe the IMEX formulation for Set2C. The linear operator for Set2C is taken as

$$L_{2C} = -\begin{pmatrix} \nabla \cdot \mathbf{U} \\ \nabla P_L' + \rho' \nabla \phi \\ \nabla \cdot \left( \dfrac{\Theta_0}{\rho_0} \mathbf{U} \right) \end{pmatrix} \tag{3.8}$$

with the linearized pressure as

$$P_L' = \frac{\gamma P_0}{\Theta_0} \Theta' \tag{3.9}$$

which is derived from a first order Taylor series approximation. It can be seen that the linear operator contains the fast moving acoustic ($P_L'$) and gravity wave ($\phi$) terms. Applying the semi-implicit formulation (Eq. (3.5)) to Set2C and letting $h_{2C} = \dfrac{\Theta_0}{\rho_0}$, yields the *No-Schur* form of the time-discretized equations

$$\rho_{tt} = \widehat{\rho} - \alpha \nabla \cdot \mathbf{U}_{tt} \tag{3.10a}$$

$$\mathbf{U}_{tt} = \widehat{\mathbf{U}} - \alpha \left( \nabla P_{tt} + \rho_{tt} \nabla \phi \right) \tag{3.10b}$$

$$\Theta_{tt} = \widehat{\Theta} - \alpha \nabla \cdot \left( h_{2C} \mathbf{U}_{tt} \right) \tag{3.10c}$$

$$P_{tt} = \frac{\gamma P_0}{\Theta_0} \Theta_{tt}. \tag{3.10d}$$

Since the linearization is performed about time-independent reference fields ($\rho_0 \neq \rho_0(t),\ \Theta_0 \neq \Theta_0(t),\ E_0 \neq E_0(t)$), the Jacobian is constant in time and only needs to be constructed once. To derive the Schur form of Eq. (3.10), we begin by substituting Eq. (3.10c) into Eq. (3.10d) to obtain

$$P_{tt} = \frac{\gamma P_0}{\Theta_0} \left( \widehat{\Theta} - \alpha \nabla \cdot \left[ h_{2C} \mathbf{U}_{tt} \right] \right). \tag{3.11}$$

Multiplying Eq. (3.10a) by $h_{2C}$ and subtracting from Eq. (3.10c) yields

$$\begin{aligned} \Theta_{tt} - h_{2C} \rho_{tt} &= \widehat{\Theta} - h_{2C} \widehat{\rho} - \alpha \left( \nabla \cdot \left( h_{2C} \mathbf{U}_{tt} \right) - h_{2C} \nabla \cdot \mathbf{U}_{tt} \right) \\ \Theta_{tt} - h_{2C} \rho_{tt} &= \widehat{\Theta} - h_{2C} \widehat{\rho} - \alpha \left( \mathbf{U}_{tt} \cdot \nabla h_{2C} \right) \end{aligned} \tag{3.12}$$

where we have made the assumption that the product rule is satisfied (although for element-based Galerkin methods this is not the case when using inexact integration - see, e.g., [1]); we will see in the results section that this assumption is reasonable. Substituting Eq. (3.10d) into Eq. (3.12) yields

$$\rho_{tt} = \widehat{\rho} - \frac{\widehat{\Theta}}{h_{2C}} + \frac{\Theta_0 P_{tt}}{h_{2C} \gamma P_0} + \frac{\alpha}{h_{2C}} \left( \mathbf{U}_{tt} \cdot \nabla h_{2C} \right). \tag{3.13}$$

Letting $G_{2C} = \dfrac{\nabla \phi}{h_{2C}}$, substituting Eq. (3.13) into Eq. (3.10b) and solving for $\mathbf{U}_{tt}$ yields

$$\mathbf{U}_{tt} = \mathbf{A}^{-1} \left( \widehat{\mathbf{U}} - \alpha \left[ \nabla P_{tt} + \left( \widehat{\rho} \nabla \phi - G_{2C} \widehat{\Theta} + G_{2C} \frac{\Theta_0 P_{tt}}{\gamma P_0} \right) \right] \right) \tag{3.14}$$

where $\mathbf{A} = 1 + \alpha^2 G_{2C} \nabla h_{2C}$. Substituting Eq. (3.14) into Eq. (3.11) and rearranging yields the Schur form for Set2C

$$P_{tt} - \frac{\gamma P_0 \alpha^2}{\Theta_0} \nabla \cdot \left( h_{2C} \mathcal{L}_{2C} \right) = \frac{\gamma P_0}{\Theta_0} \left( \widehat{\Theta} - \alpha \nabla \cdot \left[ h_{2C} \mathcal{R}_{2C} \right] \right) \tag{3.15}$$

where $\mathcal{L}_{2C} = \mathbf{A}^{-1} \left( \nabla P_{tt} + G_{2C} \frac{\Theta_0 P_{tt}}{\gamma P_0} \right)$ and $\mathcal{R}_{2C} = \mathbf{A}^{-1} \left( \widehat{\mathbf{U}} - \alpha \left[ \widehat{\rho} \nabla \phi - G_{2C} \widehat{\Theta} \right] \right)$.

This reduces the system of equations from five prognostic variables to a single (Helmholtz-like) equation for pressure.

Once the pressure is computed, the velocities can be extracted using Eq. (3.14), followed by the density Eq. (3.13), and then the potential temperature Eq. (3.12).

#### 3.2.2 1D-IMEX Formulation

The formulation in Eq. (3.15), which treats the acoustics terms in all directions, can be further reduced to implicitly treat only the vertical/radial direction. This uni-directional formulation will be referred to as 1D-IMEX and is obtained by replacing all three-dimensional (3D) differential operators $\nabla = \left( \frac{\partial}{\partial x}, \frac{\partial}{\partial y}, \frac{\partial}{\partial z} \right)$ with radial derivatives $\nabla_r = \frac{\partial}{\partial r}$ that are one-dimensional (1D). Letting $\mathbf{A}_r = 1 + \alpha^2 G_{2C} \frac{\partial h_{2C}}{\partial r}$, the 1D-IMEX formulation is given as

$$P_{tt} - \frac{\gamma P_0 \alpha^2}{\Theta_0} \nabla_r \cdot \left( h_{2C} \mathcal{L}^r_{2C} \right) = \frac{\gamma P_0}{\Theta_0} \left( \widehat{\Theta} - \alpha \nabla_r \cdot \left[ h_{2C} \mathcal{R}^r_{2C} \right] \right) \tag{3.16}$$

where $\mathcal{L}^r_{2C} = \mathbf{A}_r{}^{-1} \left( \frac{\partial P_{tt}}{\partial r} + G_{2C} \frac{\Theta_0 P_{tt}}{\gamma P_0} \right)$, and $\mathcal{R}^r_{2C} = \mathbf{A}_r{}^{-1} \left( \widehat{\mathbf{U}_r} - \alpha \left[ \widehat{\rho} \nabla \phi - G_{2C} \widehat{\Theta} \right] \right)$. It should be noted that the Schur complement is valid when $\mathbf{A}$ is nonsingular, a condition that is satisfied for a stable stratified reference atmosphere used in this work (see [23]).

### 3.3 IMplicit-EXplicit Formulation: Set3C

#### 3.3.1 3D-IMEX Formulation

We now present the IMEX formulation for Set3C. The linear operator for Set3C is taken as

$$L_{3C} = - \begin{pmatrix} \nabla \cdot \mathbf{U} \\ \nabla P'_L + \rho' \nabla \phi \\ \nabla \cdot \left( \frac{(E_0 + P_0)}{\rho_0} \mathbf{U} \right) \end{pmatrix} \tag{3.17}$$

with linearized pressure

$$P'_L = (\gamma - 1) \left( E' - \rho' \phi \right). \tag{3.18}$$

Then applying the semi-implicit formulation (Eq. (3.5)) to Set3C and letting $h_{3C} = \frac{(E_0 + P_0)}{\rho_0}$ be the reference enthalpy, yields the *No-Schur* form of the time-discretized equations

$$\rho_{tt} = \widehat{\rho} - \alpha \nabla \cdot \mathbf{U}_{tt} \tag{3.19a}$$

$$\mathbf{U}_{tt} = \widehat{\mathbf{U}} - \alpha \left( \nabla P_{tt} + \rho_{tt} \nabla \phi \right) \tag{3.19b}$$

$$E_{tt} = \widehat{E} - \alpha \nabla \cdot \left( h_{3C} \mathbf{U}_{tt} \right) \tag{3.19c}$$

$$P_{tt} = (\gamma - 1) \left( E_{tt} - \rho_{tt} \phi \right). \tag{3.19d}$$

To derive the Schur form of Eq. (3.19), we begin by substituting Eq. (3.19a) and Eq. (3.19c) into Eq. (3.19d), and letting $F_{3C} = h_{3C} - \phi$ to obtain

$$P_{tt} = (\gamma - 1)\left(\widehat{E} - \phi\widehat{\rho}\right) - \alpha(\gamma - 1)\left(F_{3C}\nabla \cdot \mathbf{U}_{tt} + \nabla h_{3C} \cdot \mathbf{U}_{tt}\right). \tag{3.20}$$

Multiplying Eq. (3.19a) by $h_{3C}$ and subtracting from Eq. (3.19c) yields

$$E_{tt} - h_{3C}\rho_{tt} = \widehat{E} - h_{3C}\widehat{\rho} - \alpha\nabla h_{3C} \cdot \mathbf{U}_{tt}. \tag{3.21}$$

Substituting Eq. (3.21) into Eq. (3.19d) yields

$$P_{tt} = (\gamma - 1)\left(\widehat{E} + h_{3C}\rho_{tt} - h_{3C}\widehat{\rho} - \alpha\nabla h_{3C} \cdot \mathbf{U}_{tt} - \rho_{tt}\phi\right) \tag{3.22}$$

which can be solved for density $\rho_{tt}$ as

$$\rho_{tt} = \frac{1}{F_{3C}}\left(\frac{1}{\gamma - 1}P_{tt} + \alpha\nabla h_{3C} \cdot \mathbf{U}_{tt} + h_{3C}\widehat{\rho} - \widehat{E}\right). \tag{3.23}$$

Letting $G_{3C} = \dfrac{\nabla\phi}{F_{3C}}$, substituting Eq. (3.23) into Eq. (3.19b), and solving for $\mathbf{U}_{tt}$ yields

$$\mathbf{U}_{tt} = \mathbf{A}^{-1}\left(\widehat{\mathbf{U}} - \alpha G_{3C}\left(h_{3C}\widehat{\rho} - \widehat{E}\right) - \alpha\nabla P_{tt} - \frac{\alpha G_{3C}}{(\gamma - 1)}P_{tt}\right) \tag{3.24}$$

where $\mathbf{A} = 1 + \alpha^2 G_{3C}\nabla h_{3C}$. Then, substituting Eq. (3.24) into Eq. (3.20) and rearranging yields the Schur form for Set3C

$$\begin{aligned} P_{tt} &- \alpha^2(\gamma - 1)F_{3C}\nabla \cdot \mathcal{L}_{3C} - \alpha^2(\gamma - 1)\nabla h_{3C} \cdot \mathcal{L}_{3C} \\ &= (\gamma - 1)\left(\widehat{E} - \phi\widehat{\rho}\right) - \alpha(\gamma - 1)F_{3C}\nabla \cdot \mathcal{R}_{3C} - \alpha(\gamma - 1)\nabla h_{3C} \cdot \mathcal{R}_{3C} \end{aligned} \tag{3.25}$$

where $\mathcal{L}_{3C} = \mathbf{A}^{-1}\left(\nabla P_{tt} + \dfrac{G_{3C}}{(\gamma - 1)}P_{tt}\right)$ and $\mathcal{R}_{3C} = \mathbf{A}^{-1}\left(\widehat{\mathbf{U}} - \alpha G_{3C}\left(h_{3C}\widehat{\rho} - \widehat{E}\right)\right)$. The Schur form can be written compactly as Eq. (3.25),

$$\begin{aligned} P_{tt} &- \alpha^2(\gamma - 1)\left(\nabla \cdot (h_{3C}\mathcal{L}_{3C}) - \phi\nabla \cdot \mathcal{L}_{3C}\right) \\ &= (\gamma - 1)\left(\widehat{E} - \phi\widehat{\rho} - \alpha\left(\nabla \cdot (h_{3C}\mathcal{R}_{3C}) - \phi\nabla \cdot \mathcal{R}_{3C}\right)\right) \end{aligned} \tag{3.26}$$

and is the form that will be used for deriving the discrete forms.

### 3.3.2 1D-IMEX Formulation

As was done for Set2C (Eq. (3.16)), the 1D-IMEX formulation for Set3C can be obtained by substituting $\nabla = \left(\dfrac{\partial}{\partial x}, \dfrac{\partial}{\partial y}, \dfrac{\partial}{\partial z}\right)$ in Eq. (3.26) with radial derivatives $\nabla_r = \dfrac{\partial}{\partial r}$. Letting $\mathbf{A}_r = 1 + \alpha^2 G_{3C}\dfrac{\partial h_{3C}}{\partial r}$, the 1D-IMEX formulation for Set3C is written as

$$\begin{aligned} P_{tt} &- \alpha^2(\gamma - 1)\left(\nabla_r \cdot (h_{3C}\mathcal{L}^r_{3C}) - \phi\nabla_r \cdot \mathcal{L}^r_{3C}\right) \\ &= (\gamma - 1)\left(\widehat{E} - \phi\widehat{\rho} - \alpha\left(\nabla_r \cdot (h_{3C}\mathcal{R}^r_{3C}) - \phi\nabla_r \cdot \mathcal{R}^r_{3C}\right)\right) \end{aligned} \tag{3.27}$$

where $\mathcal{L}^r_{3C} = {\mathbf{A}_r}^{-1}\left(\dfrac{\partial P_{tt}}{\partial r} + \dfrac{G_{3C}}{(\gamma - 1)}P_{tt}\right)$, and $\mathcal{R}^r_{3C} = {\mathbf{A}_r}^{-1}\left(\widehat{\mathbf{U}_r} - \alpha G_{3C}\left(h_{3C}\widehat{\rho} - \widehat{E}\right)\right)$.

# 4 Numerical Implementation

## 4.1 Spatial Discretization

We now turn to the description of the discontinuous Galerkin discretization ( [1]) of the Schur forms for Set2C and Set3C. Multiplying Eq. (3.15) and Eq. (3.26) by a test function $\psi$ and integrating over an element $\Omega_e$

$$\int_{\Omega_e} \psi P_{tt}\, d\Omega_e - \alpha^2 \int_{\Omega_e} \frac{\gamma P_0}{\Theta_0} \nabla \cdot (h_{2C} \mathcal{L}_{2C})\, d\Omega_e = \int_{\Omega_e} \psi \frac{\gamma P_0}{\Theta_0} \widehat{\Theta}\, d\Omega_e - \alpha \int_{\Omega_e} \psi \frac{\gamma P_0}{\Theta_0} \nabla \cdot (h_{2C} \mathcal{R}_{2C})\, d\Omega_e \tag{4.1}$$

$$\begin{aligned} \int_{\Omega_e} \psi P_{tt}\, d\Omega_e - \alpha^2 (\gamma - 1) \int_{\Omega_e} \psi \left(\nabla \cdot (h_{3C} \mathcal{L}_{3C}) - \phi \nabla \cdot \mathcal{L}_{3C}\right)\, d\Omega_e = (\gamma - 1) \int_{\Omega_e} \psi \left(\widehat{E} - \phi \widehat{\rho}\right)\, d\Omega_e \\ - \alpha (\gamma - 1) \int_{\Omega_e} \psi \left(\nabla \cdot (h_{3C} \mathcal{R}_{3C}) - \phi \nabla \cdot \mathcal{R}_{3C}\right)\, d\Omega_e, \end{aligned} \tag{4.2}$$

we obtain the variational forms of the Schur complements for Set2C and Set3C, respectively. Without loss of generality, the $\mathcal{L}$ and $\mathcal{R}$ operators for both equation sets can be written as

$$\begin{aligned} \mathcal{L} &= \mathbf{A}^{-1} \left(\nabla P_{tt} + \boldsymbol{a} P_{tt}\right) \\ \mathcal{R} &= \mathbf{A}^{-1} \left(\widehat{\mathbf{U}} - \alpha \left(\boldsymbol{b} \widehat{\rho} - \boldsymbol{c} \widehat{T}\right)\right) \end{aligned} \tag{4.3}$$

where $\boldsymbol{a} = G_{2C} \frac{\Theta_0}{\gamma P_0}$, $\boldsymbol{b} = \nabla \phi$, $\boldsymbol{c} = G_{2C}$ and the thermodynamic variable $\widehat{T} = \widehat{\Theta}$ for Set2C, and $\boldsymbol{a} = \frac{G_{3C}}{\gamma - 1}$, $\boldsymbol{b} = G_{3C} h_{3C}$, $\boldsymbol{c} = 1$ and $\widehat{T} = \widehat{E}$ for Set3C.

It can be seen from Eq. (4.1) and Eq. (4.2) that the differential operators are similar, therefore, we will show the discretization for Set3C with $\mathcal{F} \in \{\mathcal{L}_{3C}, \mathcal{R}_{3C}\}$ and $h_C = h_{3C}$ but discretizations can be obtained for Set2C with $\mathcal{F} \in \{\mathcal{L}_{2C}, \mathcal{R}_{2C}\}$, $\phi = 0$ and $h_C = h_{2C}$. Similarly, the 1D-IMEX discretization can be obtained by retaining only the radial differential operators. It can be seen that $\nabla \cdot \mathcal{L}$ results in both first and second order operators while $\nabla \cdot \mathcal{R}$ results in only first order operators. The second order operators are discretized using the local discontinuous Galerkin (LDG) method [24].

We begin by first discretizing the first-order operators in the Schur form

$$\int_{\Omega_e} \psi \left(\omega \nabla \cdot (h_C \mathcal{F}) - \phi \nabla \cdot \mathcal{F}\right)\, d\Omega_e = \int_{\Omega_e} \psi \omega \nabla \cdot (h_C \mathcal{F})\, d\Omega_e - \int_{\Omega_e} \psi \phi \nabla \cdot \mathcal{F}\, d\Omega_e \tag{4.4}$$

where $\omega = \frac{\gamma P_0}{\Theta_0}$ for Set2C, $\omega = 1$ for Set3C and $\psi$ are the multi-dimensional Lagrange polynomials defined as a tensor product of 1D Lagrange polynomials supported at Legendre-Gauss-Lobotto (LGL) points. The integrals are approximated by means of inexact integration using the LGL quadrature (for more details see, e.g., [1]). Applying integration by parts to the first term on the right-hand side yields

$$\begin{aligned} \int_{\Omega_e} \psi \omega \nabla \cdot (h_C \mathcal{F})\, d\Omega_e &= \int_{\Omega_e} \nabla \cdot (\psi \omega h_C \mathcal{F})\, d\Omega_e - \int_{\Omega_e} h_C \mathcal{F} \cdot \nabla(\psi \omega)\, d\Omega_e \\ &= \int_{\Gamma_e} \widetilde{\psi \omega h_C \mathcal{F}} \cdot \hat{\mathbf{n}}\, d\Gamma_e - \int_{\Omega_e} \omega h_C \nabla \psi \cdot \mathcal{F}\, d\Omega_e - \int_{\Omega_e} \psi h_C \nabla \omega \cdot \mathcal{F}\, d\Omega_e \end{aligned} \tag{4.5}$$

where $\widetilde{\psi \omega h_C \mathcal{F}}$ is some numerical flux. Due to the symmetry of the LGL nodes, this simplifies to $\widetilde{\psi \omega h_C \mathcal{F}} = \psi \widetilde{\omega h_C \mathcal{F}}$, i.e., the basis function satisfies $\psi = \widetilde{\psi}$. The second term on the right-hand side in Eq. (4.4) can be written as

$$\begin{aligned} \int_{\Omega_e} \psi \phi \nabla \cdot \mathcal{F}\, d\Omega_e &= \int_{\Omega_e} \nabla \cdot (\psi \phi \mathcal{F})\, d\Omega_e - \int_{\Omega_e} \mathcal{F} \cdot \nabla(\psi \phi)\, d\Omega_e \\ &= \int_{\Gamma_e} \psi \widetilde{\phi \mathcal{F}} \cdot \hat{\mathbf{n}}\, d\Gamma_e - \int_{\Omega_e} \phi \nabla \psi \cdot \mathcal{F}\, d\Omega_e - \int_{\Omega_e} \psi \nabla \phi \cdot \mathcal{F}\, d\Omega_e. \end{aligned} \tag{4.6}$$

In flow-in-a-box mode, since $\phi = gz$ then $\nabla\phi = g\widehat{k}$ where $\widehat{k}$ is the unit vector along the z direction. Substituting Eq. (4.5) and Eq. (4.6) into Eq. (4.4), and rearranging yields the weak form of the differential operator

$$\begin{aligned}\int_{\Omega_e} \psi\left(\omega\nabla\cdot\left(h_C\mathcal{F}\right)-\phi\nabla\cdot\mathcal{F}\right)\,d\Omega_e &= \int_{\Gamma_e}\psi\left(\widetilde{\omega h_C\mathcal{F}}-\widetilde{\phi\mathcal{F}}\right)\cdot\hat{\mathbf{n}}\,d\Gamma_e-\int_{\Omega_e}\left(\omega h_C-\phi\right)\nabla\psi\cdot\mathcal{F}\,d\Omega_e\\ &\quad-\int_{\Omega_e}\psi\left(h_C\nabla\omega-\nabla\phi\right)\cdot\mathcal{F}\,d\Omega_e.\end{aligned} \tag{4.7}$$

Applying integration by parts to the second term on the right-hand side of Eq. (4.7) yields

$$\begin{aligned}\int_{\Omega_e}\left(\omega h_C-\phi\right)\nabla\psi\cdot\mathcal{F}\,d\Omega_e &= \int_{\Omega_e}\nabla\cdot\left(\psi\left(\omega h_C-\phi\right)\mathcal{F}\right)d\Omega_e-\int_{\Omega_e}\psi\nabla\cdot\left(\left(\omega h_C-\phi\right)\mathcal{F}\right)d\Omega_e\\ &= \int_{\Gamma_e}\psi\left(\omega^e h_C^e-\phi^e\right)\mathcal{F}^e\cdot\hat{\mathbf{n}}\,d\Gamma_e-\int_{\Omega_e}\psi\left(\left(\omega h_C-\phi\right)\nabla\cdot\mathcal{F}+\nabla(\omega h_C-\phi)\cdot\mathcal{F}\right)\,d\Omega_e,\end{aligned} \tag{4.8}$$

where we have added the superscript $e$ to the surface integral only (first term on the right-hand side) to denote a local element-wise quantity. Substituting Eq. (4.8) into Eq. (4.7), we obtain the strong form

$$\begin{aligned}\int_{\Omega_e}\psi\left(\nabla\cdot\left(\omega h_C\mathcal{F}\right)-\phi\nabla\cdot\mathcal{F}\right)\,d\Omega_e &= \int_{\Gamma_e}\psi\left(\left(\widetilde{\omega h_C\mathcal{F}}-\widetilde{\phi\mathcal{F}}\right)-\left(\omega^e h_C^e-\phi^e\right)\mathcal{F}^e\right)\cdot\hat{\mathbf{n}}\,d\Gamma_e\\ &\quad+\int_{\Omega_e}\psi\left(\omega h_C-\phi\right)\nabla\cdot\mathcal{F}\,d\Omega_e+\int_{\Omega_e}\psi\omega\nabla h_C\cdot\mathcal{F}\,d\Omega_e\end{aligned} \tag{4.9}$$

where it should be understood that terms with no designation (either a superscript $e$ or tilde) strictly represent local element-wise quantities. If we assume the reference fields to be continuous then $\left(\widetilde{\omega h_C\mathcal{F}}-\widetilde{\phi\mathcal{F}}\right)=\left(\omega h_C-\phi\right)\widetilde{\mathcal{F}}$. Furthermore, if $\mathcal{F}=\mathcal{R}$ then

$$\widetilde{\mathcal{R}}=\mathbf{A}^{-1}\left(\widetilde{\widehat{\mathbf{U}}}-\alpha\left(\boldsymbol{b}\,\widetilde{\widehat{\rho}}-\boldsymbol{c}\,\widetilde{\widehat{T}}\right)\right)$$

since $\mathbf{A}$, $\boldsymbol{a}$, $\boldsymbol{b}$ and $\boldsymbol{c}$ are also continuous (as they are dependent on the reference fields). The case where $\mathcal{F}=\mathcal{L}$ results in both first and second order operators. Expanding $\mathcal{L}=\mathbf{A}^{-1}\nabla P_{tt}+\mathbf{A}^{-1}\boldsymbol{a}P_{tt}$, the second term on the right-hand side can be treated by letting $\mathcal{F}=\mathbf{A}^{-1}\boldsymbol{a}P_{tt}$ in Eq. (4.7) and Eq. (4.9) (and setting the coefficients appropriately), but the first term requires computation of the variational form of the gradient consistent with the DG scheme. This is done by multiplying $\nabla P_{tt}$ by a tensor test function $\Psi=\psi\mathbf{I}$ and invoking the divergence theorem resulting in the variational weak form of the gradient

$$\begin{aligned}\int_{\Omega_e}\Psi\cdot\nabla P_{tt}\,d\Omega_e &= \int_{\Omega_e}\nabla\cdot\left(\Psi P_{tt}\right)d\Omega_e-\int_{\Omega_e}P_{tt}\nabla\cdot\Psi\,d\Omega_e\\ &= \int_{\Gamma_e}\widetilde{P_{tt}\Psi}\cdot\hat{\mathbf{n}}\,d\Gamma_e-\int_{\Omega_e}P_{tt}\nabla\cdot\Psi\,d\Omega_e\end{aligned} \tag{4.10}$$

where $\widetilde{P_{tt}\Psi}=\widetilde{P_{tt}}\Psi$, since $\widetilde{\psi}=\psi$. Applying integration by parts to the second term on the right hand side of Eq. (4.10) yields

$$\begin{aligned}\int_{\Omega_e}P_{tt}\nabla\cdot\Psi\,d\Omega_e &= \int_{\Omega_e}\nabla\cdot\left(P_{tt}\Psi\right)d\Omega_e-\int_{\Omega_e}\nabla P_{tt}\cdot\Psi\,d\Omega_e\\ &= \int_{\Gamma_e}P_{tt}^e\Psi\cdot\hat{\mathbf{n}}\,d\Gamma_e-\int_{\Omega_e}\nabla P_{tt}\cdot\Psi\,d\Omega_e.\end{aligned} \tag{4.11}$$

Substituting Eq. (4.11) into Eq. (4.10) yields the variational strong form for the gradient

$$
\begin{aligned}
\int_{\Omega_e} \Psi \cdot \nabla P_{tt}\, d\Omega_e &= \int_{\Omega_e} \nabla \cdot (\Psi P_{tt})\, d\Omega_e - \int_{\Omega_e} P_{tt} \nabla \cdot \Psi\, d\Omega_e \\
&= \int_{\Gamma_e} \left(\widetilde{P_{tt}} - P_{tt}^e\right) \Psi \cdot \hat{\mathbf{n}}\, d\Gamma_e + \int_{\Omega_e} \nabla P_{tt} \cdot \Psi\, d\Omega_e.
\end{aligned} \tag{4.12}
$$

Then the weak and strong forms for $\mathcal{L}$ are

$$
\begin{aligned}
\mathcal{L}_{weak} &= \mathbf{A}^{-1}\left[\int_{\Gamma_e} \widetilde{P_{tt}} \Psi \cdot \hat{\mathbf{n}}\, d\Gamma_e - \int_{\Omega_e} P_{tt} \nabla \cdot \Psi\, d\Omega_e + \boldsymbol{a} P_{tt}\right] \\
&= \mathbf{A}^{-1}\left[\boldsymbol{\mathcal{C}}^{(e,k)} - \widetilde{\boldsymbol{D}}^{(e)} + \boldsymbol{a}\mathbf{I}\right] P_{tt}
\end{aligned} \tag{4.13}
$$

and

$$
\begin{aligned}
\mathcal{L}_{strong} &= \mathbf{A}^{-1}\left[\int_{\Gamma_e} \left(\widetilde{P_{tt}} - P_{tt}^e\right) \Psi \cdot \hat{\mathbf{n}}\, d\Gamma_e + \int_{\Omega_e} \nabla P_{tt} \cdot \Psi\, d\Omega_e + \boldsymbol{a} P_{tt}\right] \\
&= \mathbf{A}^{-1}\left[\boldsymbol{\Delta}^{(e,k)} + \boldsymbol{D}^{(e)} + \boldsymbol{a}\mathbf{I}\right] P_{tt}
\end{aligned} \tag{4.14}
$$

respectively, where $\mathbf{I}$ is the identity matrix, $\boldsymbol{\mathcal{C}}^{(e,k)}$ and $\boldsymbol{\Delta}^{(e,k)}$ are the centered and difference matrices, and $\boldsymbol{D}^{(e)}$ and $\widetilde{\boldsymbol{D}}^{(e)}$ are the strong and weak form differentiation matrices, which are defined in Appendix C with more details found in [ [1], Ch. 21]. The matrices with superscript $(e, k)$ are global matrices that required communication between the element $e$ and its neighbor $k$, whereas the matrices with superscript $(e)$ are local matrices that do not require communication with neighboring elements. Since the fully discrete representation of the Schur forms are unwieldy, they are presented in Appendix C. Unless explicitly stated the use of the strong form, the results and analysis of various test cases are obtained using the weak form.

### 4.2 Conditions on Numerical Flux for Consistent IMEX Splitting

Here we show the conditions required for the IMEX splitting to be consistent on the continuous and discrete levels. We only show this for Set3C, but the same conditions apply for Set2C. Consider a balance law (Set3C) written in vector form defined by Eq. (2.1) and Eq. (2.5), where, for clarity, we denote the total flux as $\mathbf{F}_T = \mathbf{F}$ with components $(F^\rho, \mathbf{F}^\mathbf{U}, F^E)$ that represent the mass, momentum and total energy flux. The coloring and subscripts $T$, $L$, and $NL$ denote total, linear, and nonlinear components, respectively.

Here onwards, for convenience, the source terms $\mathbf{S}$ due to gravity and Coriolis are neglected; thereby effectively replacing the original balance law with a conservation law with no loss in generality (because the stiffness in the equations are derived from the terms in the flux terms). Casting Eq. (2.1) into the variational weak form (and invoking the divergence theorem) yields

$$
\int_{\Omega_e} \psi \frac{\partial \mathbf{q}}{\partial t}\, d\Omega_e + \int_{\Gamma_e} \widetilde{\psi \mathbf{F}_T} \cdot \hat{n}\, d\Gamma_e - \int_{\Omega_e} \nabla\psi \cdot \mathbf{F}_T\, d\Omega_e = 0 \tag{4.15}
$$

where $\widetilde{Q}$ represents the numerical flux for the quantity $Q$. Due to the symmetry of the Legendre-Gauss-Lobotto (LGL) points, tensor product representation of the Lagrange polynomials, and conforming elements, $\widetilde{\psi \mathbf{F}_T} = \psi\widetilde{\mathbf{F}_T}$. Using the Rusanov flux, Eq. (4.15) then becomes

$$
\int_{\Omega_e} \psi \frac{\partial \mathbf{q}}{\partial t}\, d\Omega_e + \int_{\Gamma_e} \psi \left(\{\{\mathbf{F}_T\}\} - \frac{\lambda_T}{2} [\![\mathbf{q}]\!] \hat{n}\right) \cdot \hat{n}\, d\Gamma_e - \int_{\Omega_e} \nabla\psi \cdot \mathbf{F}_T\, d\Omega_e = 0 \tag{4.16}
$$

where $\lambda_T = (u + a)^\star$ is the maximum eigenvalue (wave speed) of the of total flux Jacobian $\left(\frac{\partial \mathbf{F}_T}{\partial \mathbf{q}}\right)$, $\star$ denotes the maximum of the left and right states (i.e., $A^\star = \max(\mid A^L \mid, \mid A^R \mid)$), and $a = \sqrt{\frac{\gamma P}{\rho}}$ is the speed of sound.

The flux vector $\mathbf{F}_T := \mathbf{F}$ in Eq. (2.5) on the continuous level can be decomposed into its linear $\mathbf{F}_L$ and nonlinear $\mathbf{F}_{NL} = \mathbf{F}_T - \mathbf{F}_L$ components as

$$\mathbf{F}_L = \begin{Bmatrix} \mathbf{U} \\ P_L \mathbf{I} \\ h_0 \mathbf{U} \end{Bmatrix}, \quad \mathbf{F}_{NL} = \begin{Bmatrix} 0 \\ \frac{\mathbf{U} \otimes \mathbf{U}}{\rho} \\ \left( \frac{E+P}{\rho} - h_0 \right) \mathbf{U} \end{Bmatrix}. \tag{4.17}$$

where $h_0 = \frac{E_0 + P_0}{\rho_0}$ is the reference enthalpy and $P_L = (\gamma - 1)(E - \rho\phi)$. Using such a (simplified) decomposition and rearranging terms, Eq. (4.16) can also be decomposed into its linear (blue) and nonlinear components (red) on the discrete level as

$$\begin{aligned} \int_{\Omega_e} \psi \frac{\partial \mathbf{q}}{\partial t} \, d\Omega_e + \underbrace{\int_{\Gamma_e} \psi \left( \{\{\mathbf{F}_{NL}\}\} - \frac{\lambda_T}{2} [\![\mathbf{q}]\!] \widehat{\mathbf{n}} \right) \cdot \widehat{\mathbf{n}} \, d\Gamma_e - \int_{\Omega_e} \nabla\psi \cdot \mathbf{F}_{NL} \, d\Omega_e}_{\int_{\Omega_e} \psi \nabla \cdot \mathbf{F}_{NL} \, d\Omega_e} \\ + \underbrace{\int_{\Gamma_e} \psi \{\{\mathbf{F}_L\}\} \cdot \widehat{\mathbf{n}} \, d\Gamma_e - \int_{\Omega_e} \nabla\psi \cdot \mathbf{F}_L \, d\Omega_e}_{\int_{\Omega_e} \psi \nabla \cdot \mathbf{F}_L \, d\Omega_e} = 0 \end{aligned} \tag{4.18}$$

where the linear (blue) term represents the variational weak form of the linear flux terms using centered fluxes and the nonlinear (red) term represents the variational weak form of the nonlinear flux terms using the Rusanov flux. The linear pressure term appears nonlinearly in the wave speed $\lambda_T = \lambda(u, a) = \lambda(u, P, \rho)$ in the treatment of the nonlinear fluxes and therefore, the two terms are not separable.

Without loss of generality, we can write $\lambda_T = \lambda_L + \lambda_{NL}$ to obtain a form where each component can be treated using the Rusanov flux using the appropriate wave speeds,

$$\begin{aligned} \int_{\Omega_e} \psi \frac{\partial \mathbf{q}}{\partial t} \, d\Omega_e + \underbrace{\int_{\Gamma_e} \psi \left( \{\{\mathbf{F}_{NL}\}\} - \frac{\lambda_{NL}}{2} [\![\mathbf{q}]\!] \widehat{\mathbf{n}} \right) \cdot \widehat{\mathbf{n}} \, d\Gamma_e - \int_{\Omega_e} \nabla\psi \cdot \mathbf{F}_{NL} \, d\Omega_e}_{\int_{\Omega_e} \psi \nabla \cdot \mathbf{F}_{NL} \, d\Omega_e} \\ + \underbrace{\int_{\Gamma_e} \psi \left( \{\{\mathbf{F}_L\}\} - \frac{\lambda_L}{2} [\![\mathbf{q}]\!] \widehat{\mathbf{n}} \right) \cdot \widehat{\mathbf{n}} \, d\Gamma_e - \int_{\Omega_e} \nabla\psi \cdot \mathbf{F}_L \, d\Omega_e}_{\int_{\Omega_e} \psi \nabla \cdot \mathbf{F}_L \, d\Omega_e} = 0 \end{aligned} \tag{4.19}$$

where $\lambda_L$ and $\lambda_{NL}$ are now the maximum eigenvalues of the Jacobians of the linear flux $\mathbf{F}_L$ and nonlinear flux $\mathbf{F}_{NL}$, respectively.

In one dimension, the linear and total flux terms are given as follows

$$\mathbf{F}_L = \begin{Bmatrix} U \\ P_L \\ h_0 U \end{Bmatrix}, \quad \mathbf{F}_T = \begin{Bmatrix} U \\ \frac{U^2}{\rho} \\ \left( \frac{E+P}{\rho} \right) U \end{Bmatrix} \tag{4.20}$$

where the buoyancy term in the momentum equation is omitted because it is considered a source. Computing the Jacobian of the linear flux yields

$$J_L = \frac{\partial \mathbf{F}_L}{\partial \mathbf{q}} = \begin{pmatrix} 0 & 1 & 0 \\ -(\gamma - 1)\phi & 0 & \gamma - 1 \\ 0 & h_0 & 0 \end{pmatrix} \tag{4.21}$$

with eigenvalues $\lambda_L = (-a_0, 0, a_0)$ where $a_0 = \sqrt{\frac{\gamma P_0}{\rho_0}}$. Computing the Jacobian of the total flux $\mathbf{F}_T$ yields

$$J_T = \frac{\partial \mathbf{F}_T}{\partial \mathbf{q}} = \begin{pmatrix} 0 & 1 & 0 \\ (\gamma - 3)\frac{U^2}{2\rho^2} - (\gamma - 1)\phi & (3 - \gamma)\frac{U}{\rho} & \gamma - 1 \\ \frac{U}{\rho^2} \left( -\gamma E + (\gamma - 1)\frac{U}{\rho^2} \right) & \frac{\gamma E}{\rho} - \frac{3(\gamma-1)U^2}{2\rho^2} - (\gamma - 1)\phi & \frac{\gamma U}{\rho} \end{pmatrix} \tag{4.22}$$

with eigenvalues $\lambda_T = (u - a, u, u + a)$ where $a = \sqrt{\frac{\gamma p}{\rho}}$.

Therefore, the Jacobian ($\Delta J = \frac{\partial \mathbf{F}_{NL}}{\partial \mathbf{q}}$) of the nonlinear flux $\mathbf{F}_{NL} = \mathbf{F}_T - \mathbf{F}_L$ is a lower triangular matrix with eigenvalues $\Delta\lambda = \lambda_{NL} = (0, \gamma u, (3-\gamma)u)$. The important point here is that even though the operator $F_L$ has a linearized pressure term, the eigenvalues of $\Delta J$ only have velocity components. [2] Although this was shown for the case where the nonlinear operators are implicitly constructed as $\mathbf{F}_{NL} = \mathbf{F}_T - \mathbf{F}_L$, it also holds true when the nonlinear operators are obtained by first linearizing $\mathbf{F}_T$ then splitting the terms. This can be shown by using Eq. (3.18) and linearizing the total energy flux $F^E$ about $E_0$ and $P_0$, and further about $\rho_0$ as

$$
\begin{aligned}
\left(\frac{E+P}{\rho}\right) U &= \left(\frac{E_0 + P_0}{\rho}\right) U + \left(\frac{E' + P'}{\rho}\right) U \\
&\approx h_0 U - \frac{h_0}{\rho_0} U \rho' + \left(\frac{E' + P'}{\rho}\right) U \\
&\approx \underbrace{h_0 U}_{F_L^E} - \underbrace{\left(\frac{\gamma E' - (\gamma - 1)\rho'\phi}{\rho} - \frac{h_0}{\rho_0}\rho'\right) U}_{F_{NL}^E} .
\end{aligned}
\tag{4.23}
$$

Then nonlinear flux term can be written as

$$
\mathbf{F}_{NL} = \left\{0,\ \frac{U^2}{\rho},\ \left(\frac{\gamma E' - (\gamma - 1)\rho'\phi}{\rho} - \frac{h_0}{\rho_0}\rho'\right) U\right\}^{\mathcal{T}}
\tag{4.24}
$$

where the Jacobian of the nonlinear flux again results in a lower-triangular matrix but with eigenvalues $\lambda_{NL} = (0, u, \gamma u)$ that only depend on the velocity components. This shows that for the linear and nonlinear terms to be separable and the splitting to be consistent on both the continuous and discrete levels, the numerical flux of the linear and nonlinear operators should be computed with wave speeds $\lambda_L = a_0$ and $\lambda_{NL} = u$, respectively. For the case where $\mathbf{F}_{NL}$ is computed implicitly as $\mathbf{F}_{NL} = \mathbf{F}_T - \mathbf{F}_L$, the wave speed for $\mathbf{F}_T$ should be taken as $\lambda_T = \lambda_L + \lambda_{NL}$, where $\lambda_L = 0$ when the implicit term is treated with a centered flux and $\lambda_L = a_0$ when the implicit term is treated with a linearized Rusanov flux. It should be mentioned that since $u \ll a_0$ in subsonic applications, it is possible that the eigenspectra can lie outside the region of stability of the time integration scheme used to treat the nonlinear terms explicitly. To keep the IMEX formulation for Runge-Kutta methods that we presented in [14] unchanged, we construct the nonlinear operator as $\mathbf{F}_{NL} = \mathbf{F}_T - \mathbf{F}_L$. Since the nonzero eigenvalues of $\mathbf{F}_{NL}$ are all of relatively the same magnitude (i.e. $\mathcal{O}(u) = \mathcal{O}(\gamma u) = \mathcal{O}((3-\gamma)u)$) regardless of how $\mathbf{F}_{NL}$ is computed, we take $\lambda_{NL} = u$.

Table 1 presents the various flux combinations that can be used for the implicit and explicit dynamics. In the following sections, for all IMEX methods we call the *AT* flux that which uses the acoustic Rusanov flux ($\lambda_L = a_0$ ) for the linear implicit part (IM) combined with total Rusanov flux with $\lambda_T = u + a$ (and $\lambda_{NL} = u + a'$) accounting for both the advective and acoustic wave speeds for the nonlinear explicit (EX) part; in contrast, we call the *CA* flux that which uses the centered flux for the linear implicit (IM) part (i.e. $\lambda_L = 0$) combined with only the advective Rusanov flux ($\lambda_T = \lambda_{NL} = u$) for the nonlinear explicit (EX) part (see **Case 1** in Appendix A). Insignificant differences are seen between the total Rusanov (*AT* Flux) and linearized Rusanov (*AL* Flux) when constructing the total flux $\mathbf{F}_T$, the reasons for which are presented in **Case 3** in Appendix A.

The *No-Schur* system can be solved using either the *AT* flux or the *CA* flux formulation, whereas the Schur system will be solved using the *CA* flux approach only. The Schur form can also be solved using the centered flux for the linear term ($\lambda_L = 0$) and a Rusanov flux for the total/nonlinear terms (e.g., $\lambda_T = u + a_0$ as in *CL* Flux or $\lambda_T = u + a$ as in *CT* Flux); such a choice will subject the explicit time stepping to the CFL restriction due to the acoustic dynamics and negate the benefits of the IMEX methods (see **Case 2** in Appendix A). In Appendix B, we show why it is challenging to construct a Schur form when the linear term is treated with the Rusanov flux or any flux containing a jump term ($\lambda_L \neq 0$).

[2]Although we only show the results for one dimension, the two- and three-dimensional results bring us to the same conclusion.

**Table 1:** Flux combinations and wave speeds (magnitude of the *jump* term) in the Rusanov flux for the IMEX formulations

| Name | Abbrv. | $\lambda_L$ | $\lambda_{NL}$ | $\lambda_T$ |
|---|---|---|---|---|
| Acoustic-Total Rusanov | AT | $a_0$ | $u + a'$ | $u + a$ |
| Acoustic-Linearized Rusanov | AL | $a_0$ | $u$ | $u + a_0$ |
| Centered-Total Rusanov | CT | 0 | $u + a$ | $u + a$ |
| Centered-Linearized Rusanov | CL | 0 | $u + a_0$ | $u + a_0$ |
| Centered-Advective Rusanov | CA | 0 | $u$ | $u$ |

### 4.3 Numerical Stabilization

It should be noted that the complete system discretized using the *CA* flux is less dissipative than the system discretized using the *AT* flux, particularly for cases where $u \ll a$ (see **Case 1** in Appendix A).

This lack of dissipation can lead to spurious oscillations that can drive the system unstable, leading to divergence (i.e., blow-up). It was previously shown that artificial viscosity and sub-grid-scale (SGS) models can be used to stabilize the continuous Galerkin (which is also inherently non-dissipative) [25] and discontinuous Galerkin formulations [26]. For this reason, for some test problems we introduce an artificial diffusion operator (Laplacian) with uniform viscosity or SGS viscosity in the governing equations. The SGS viscosity used in this work is obtained using the Vreman model. It was shown in [25] that the Vreman model is able to stabilize the non-dissipative system while limiting the amount of artificial dissipation introduced. For brevity, the reader is referred to [27] for the complete description of the Vreman model. The second order diffusion/Laplacian operator ($\nabla \cdot (\mu \nabla \mathbf{U})$) is discretized using the weak-form local discontinuous Galerkin (LDG) method [24]. In the future, it will be interesting to explore entropy-stable methods in this context that do not require additional dissipation mechanisms (see, e.g., [21]).

### 4.4 Linear Solvers

Before the linear solvers used to solve the resulting linear systems are discussed, the domain decomposition for multi-core processing is briefly presented. As is often done in NWP or climate modeling (e.g., see [28]), the spherical manifold is partitioned using either a cubed-sphere or icosahedral mesh and extruded along the radial direction (e.g., see [29]) such that each column of elements lies completely on a single processor (Fig. 1). This approach allows for the construction of all radial operators without requiring intra-processor communication, thereby, increasing scaling efficiency [29–32]. All graph partitioning and parallel mesh decomposition is performed using the *p4est* (for 3D-IMEX) and *p6est* (for 1D-IMEX) library [33].

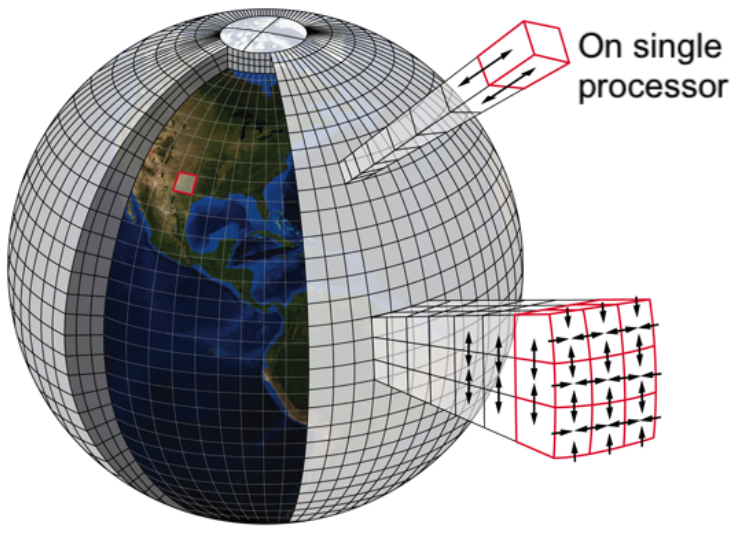

**Figure 1:** The decomposition of the computational domain for multi-core processing. [34, Ch. 2]

The choice of the formulation (*No-Schur* or Schur) can directly influence the efficiency of the solver for the linear system. Since the vertical/radial degrees-of-freedom (DOF) $N_z$ (i.e., DOF per mesh column) is much lower than the total DOF $N_{tot}$, and since the linear system in the 1D-IMEX formulation is independent of the adjacent columns, it is beneficial to use direct solvers. As stated previously, the Jacobian is constant in time and, therefore, can be factored once and stored for all subsequent linear solves. The linear system resulting from the 3D-IMEX discretization is often too large to store, too sparse to solve using direct methods, and requires intra-processor communication. For such systems, iterative Krylov methods such as GMRES, conjugate gradient (CG), or Bi-Conjugate Gradient-Stabilized (BiCGStab) are preferred. In this study we solve all 3D-IMEX systems using GMRES with a relative tolerance of $10^{-4}$, and all 1D-IMEX systems, which have a banded form, using banded LU factorization from the LAPACK library. We show in the following sections that the Schur form of the IMEX formulations results in all real, positive eigenvalues, indicating a symmetric positive-definite (SPD) system. Therefore, the Schur systems can also be solved using conjugate gradient or Cholesky factorization. It should be mentioned that for the 1D-IMEX systems, Cholesky decomposition offers negligible gains over LU factorization since the matrix is constructed and decomposed (either LU or Cholesky) only once and reused for each subsequent linear solve (i.e. forward and backward substitution).

## 5 Test Problems

This section briefly presents the four test problems used to analyze the semi-implicit discretizations. These problems are common benchmarks for verifying atmospheric models. For a more in-depth description of the test problems, readers are referred to their corresponding literature.

### 5.1 Rising Thermal Bubble

The rising thermal bubble (RTB) problem features the evolution of a warm bubble in a constant temperature ambient environment. The rising bubble deforms due to the shear stress caused by the velocity gradients. The test problem conditions [15] are similar to those of Robert [35] and features an ambient environment initially at rest and in hydrostatic balance. A warm air bubble with a potential temperature perturbation $\theta'$ is placed in the ambient atmosphere with uniform potential temperature $\theta_0$. The potential temperature perturbation distribution is defined by

$$\theta' = \begin{cases} 0, & \text{for } r > r_c \\ \dfrac{\theta_c}{2}\left(1 + \cos\left(\dfrac{\pi r}{r_c}\right)\right), & \text{for } r \leqslant r_c \end{cases} \tag{5.1}$$

where $\theta_c = 0.5\ K$, $r = \sqrt{(x - x_c)^2 + (z - z_c)^2}$ with the following constants: $\theta_0 = 300\ K$, $r_c = 250\ m$, and $(x, z) \in [0, 1000]^2\ m$ with $t \in [0, 650]\ s$ and $(x_c, z_c) = (500, 350)\ m$. The no-flux boundary condition is imposed on all four walls of the boundary.

### 5.2 Density Current

This two-dimensional test case, introduced in [36], is a standard benchmark for the development and verification of atmospheric models. It consists of a bubble of cold air descending to the ground in a neutrally stratified atmosphere. Upon hitting the lower boundary (all four boundaries are set as no-flux), the bubble develops Kelvin-Helmholtz shear instability rotors as it spreads laterally. This case is often solved by assuming a constant and uniform diffusion coefficient $\mu = 75\ m^2 s^{-1}$. The initial distribution of the potential temperature is given as

$$\theta' = \begin{cases} 0, & \text{for } r > r_c \\ \dfrac{\theta_c}{2}\left(1 + \cos\left(\dfrac{\pi r}{r_c}\right)\right), & \text{for } r \leqslant r_c \end{cases} \tag{5.2}$$

where $\theta_c = -15\ K$, $r = \sqrt{\left(\dfrac{x - x_c}{x_r}\right)^2 + \left(\dfrac{z - z_c}{z_r}\right)^2}$, and $r_c = 1$. The domain is defined as $(x, z) \in [0, 25600] \times [0, 6400]\ m$ with $t \in [0, 900]\ s$ and the center of the bubble as $(x_c, z_c) = (0, 3000)\ m$ with the size of the bubble taken to be $(x_r, z_r) = (4000, 2000)\ m$.

### 5.3 Inertia Gravity Wave

The inertia gravity wave problem features the symmetric propagation of potential temperature perturbation in a periodic channel with a constant mean flow [37, Section 3]. The initial state of the atmosphere is taken to have a constant mean flow of $\bar{u} = 20\ m/s$ with reference fields defined as

$$T_0 = 250\ K, \quad p_0(z) = p_s e^{-\delta z}, \quad \rho_0(z) = \rho_s e^{-\delta z}$$

where $p_s = 10^5\ Pa$, $\rho_s = \dfrac{p_s}{RT_0}$, and $\delta = \dfrac{g}{RT_0}$. The initial perturbation fields are given by

$$T'(x, z) = e^{\frac{\delta z}{2}} T_b(x, z), \quad p'(x, z) = 0, \quad \rho'(x, z) = e^{-\frac{\delta z}{2}} \rho_b(x, z)$$

with

$$T_b(x,z) = \Delta T \sin\left(\pi \frac{z}{h_c}\right) e^{-\left(\frac{x-x_c}{a_c}\right)^2}, \quad \rho_b(x,z) = -\rho_s \frac{T_b(x,z)}{T_0},$$

where $\Delta T = 0.001\ K$, $h_c = 10\ km$, $a_c = 5\ km$, and $x_c = 100\ km$. The domain is defined as $(x,z) \in [0,300] \times [0,10]\ km$ with $t \in [0,2500]\ s$. The left and right sides of the domain are taken to be periodic, while the top and bottom boundaries are treated as no-flux boundaries. This test problem (modified from [38]) admits an analytical solution (see [37, Section 2]).

### 5.4 Acoustic Wave on a Sphere

The acoustic wave propagation problem features an acoustic wave caused by a pressure disturbance traveling around the globe. This problem, proposed by Tomita and Satoh [39], is defined with a hydrostatically balanced state and an isothermal background potential temperature of $\theta_0 = 300\ K$. The pressure perturbation is defined as $P' = f(\lambda,\phi)g(r)$, where

$$\begin{aligned} f(\lambda,\phi) &= \begin{cases} 0, & \text{for } r > r_c \\ \dfrac{\Delta P}{2}\left(1 + \cos\left(\dfrac{\pi r}{r_c}\right)\right), & \text{for } r \leqslant r_c \end{cases} \\ g(r) &= \sin\left(\frac{n_v \pi r}{r_T}\right), \end{aligned} \tag{5.3}$$

where $\Delta P = 100\ Pa$, $n_v = 1$, $r_c = r_e/3$ is one third of the radius of the earth $r_e = 6371\ km$, and a model top of $r_T = 10\ km$. The geodesic distance $r$ is calculated as

$$r = r_e \cos^{-1}\left(\sin\phi_0 \sin\phi + \cos\phi_0 \cos\phi \cos(\lambda - \lambda_0)\right)$$

where $(\lambda_0, \phi_0)$ is the origin of the acoustic wave. The no-flux boundary condition is imposed at the top and bottom boundaries. The spatial discretization for this problem often features high aspect ratio cells where the characteristic length of the elements in the radial direction is much smaller than the length in the horizontal direction.

### 5.5 Metrics of Performance

To evaluate the performance of various semi-implicit formulations, we consider metrics that evaluate discrete conservation, efficiency, and accuracy of the methods. The methods' conservation properties can be quantified using relative mass and energy loss over simulation time $t$. The relative mass loss $\Delta\mathcal{M}$ and energy loss $\Delta\mathcal{E}$ are defined as

$$\Delta\mathcal{M} = \left|\frac{\mathcal{M}(t) - \mathcal{M}(0)}{\mathcal{M}(0)}\right| \quad \text{and} \quad \Delta\mathcal{E} = \left|\frac{\mathcal{E}(t) - \mathcal{E}(0)}{\mathcal{E}(0)}\right|$$

where $\mathcal{M}(t) = \int_\Omega \rho(\mathbf{x},t)d\Omega$ and $\mathcal{E}(t) = \int_\Omega E(\mathbf{x},t)d\Omega$. The accuracy is quantified using the mean-absolute error (MAE) relative to an analytical or reference solution $\mathbf{q}^{ref}$ as

$$\|\mathbf{q}\|_{MAE} = \frac{1}{N}\sum_i \left|\mathbf{q}_i^{ref} - \mathbf{q}_i\right|$$

where $\mathbf{q}$ is the set of prognostic variables for the equation set. Due to the lack of known analytical solutions for realistic problems in atmospheric sciences, the accuracy of the method is quantified relative to the solution obtained using an explicit scheme with a very small time-step (the small time-step is used in order to isolate the spatial error). The reference solutions are generated with a third-order, five-stage Runge-Kutta (RK35) scheme [40] with a Courant number of CN $\approx 0.0002$. We define the Courant number as

$$\text{CN} = \max\left(\frac{c\Delta t}{\Delta s}\right) \tag{5.4}$$

where $c = |u_n + a|$ is the characteristic wave speed, $u_n = \mathbf{u} \cdot \hat{\mathbf{n}}$ is the velocity in the direction $\hat{\mathbf{n}}$, $a$ is the speed of sound and $\Delta s = \sqrt{\Delta x^2 + \Delta y^2 + \Delta z^2}$ is the characteristic length. Furthermore, a separate Courant number can be specified for each of the dynamics (acoustic and advection) by considering only the appropriate wave speeds in the expression for CN. Finally, the efficiency of the formulation is studied by analyzing the time to solution and number of iterations needed to converge to the solution. For each of the different IMEX formulations, we present the condition number $\kappa(A)$ defined as

$$\kappa(A) = \frac{\sigma_{\max}}{\sigma_{\min}},$$

where $\sigma$ are the singular values of an $n \times n$ matrix $A$.

# 6 Results

## 6.1 Case 1: Rising Thermal Bubble

The rising thermal bubble problem is solved on a $30 \times 30$ mesh of fourth-order elements ($\Delta x = \Delta z = 8.3\ m$) and is evolved until $t = 650\ s$. It should be mentioned that all spatial resolutions $\Delta x$, $\Delta y$ and $\Delta z$ are defined as the average distance between the LGL nodes. The smaller dissipation in the *CA* flux IMEX discretization drives the simulation unstable whereas the added dissipation of the *AT* flux IMEX formulation keeps the model stable for the duration of the simulation. For a consistent comparison between all formulations, the Vreman SGS model, with $C_s = 0.21$, is used to stabilize all simulations of the rising thermal bubble.

The distribution of the potential temperature at $t = 650\ s$ using the various formulations of the two equation sets is shown in Fig. 2. The two equation sets result in similar distributions of the potential temperature. The 3D-IMEX and 1D-IMEX, using the same flux scheme, also result in similar solutions. There is, however, a difference when comparing results obtained using the *AT* and *CA* flux formulations. The smaller magnitude of dissipation in the *CA* flux formulation leads to a warping (Rayleigh-Taylor instabilities) at the top of the bubble, whereas the solution obtained using the *AT* flux formulation shows a bubble that is relatively intact. A similar bubble breakdown behavior is seen when using a mixed-finite element, discontinuous Galerkin formulation (see [41]).

This can be due to the *jump* term in the Rusanov flux that is absent in the *CA* flux. The results of the Schur formulation shows a similar distribution of the potential temperature as the *No-Schur* formulation for both equation sets and IMEX formulations (3D and 1D), indicating a correct derivation and implementation of the two forms. Based on the results in Fig. 2, we can conjecture that the CA flux is more appropriate in an IMEX setting since it is less dissipative (resulting in the Rayleigh-Taylor instabilities) while avoiding the production of unphysical grid imprinting (as exhibited in the AT flux).

Figure 3 shows the eigenspectra of the flux formulations, and semi-implicit discretizations. The linear system for the *No-Schur* case with the *CA* flux formulation, for both the 3D and 1D IMEX methods, has its eigenvalues along the imaginary axis. Similar spectra are seen for the IMEX formulations discretized using the continuous Galerkin method [23,42]. This is so because the continuous Galerkin formulation is inherently a spatially non-dissipative scheme. The central flux is also non-dissipative and, therefore, the eigenvalues of the IMEX formulation using continuous Galerkin and centered-flux discontinuous Galerkin methods have similar eigenspectra (e.g., see Chs. 6 and 7 in [1]).

In the case of the IMEX methods with the dissipative *AT* flux, the eigenvalues seep into the real part of the complex plane. A similar behavior was presented for discontinuous Galerkin formulations in [43] where the non-dissipative Bassi-Rebay scheme [44] resulted in only imaginary eigenvalues while the dissipative Baumann-Oden scheme [45] resulted in eigenvalues spread across the complex plane.

The eigenvalues for both the 3D and 1D Schur formulations lie entirely along the real axis - as it should for an elliptic differential operator. Although not shown, the eigenspectra of Set2C are similar to those of Set3C for all formulations. The condition number as a function of Courant number for Set2C and Set3C are shown in Figs. 3e and 3f, respectively. The *No-Schur* linear system for Set2C is much better conditioned than that of Set3C; this can be attributed to $h_{2C} \ll h_{3C}$ in Eq. (3.8) and Eq. (3.17), respectively.

Comparing the Schur and *No-Schur* formulations for Set2C, the Schur formulation results in a system with a smaller condition number for lower Courant numbers, after which the condition number grows larger than that of the *No-Schur* system, particularly for the 1D-IMEX formulation. Comparing the Schur and *No-Schur* formulations for Set3C, the linear system for the Schur form is consistently better conditioned with condition numbers five orders of magnitude smaller than that of its *No-Schur* counterpart. This is expected to significantly decrease the number of iteration required to reach convergence, thereby reducing the time to solution. The condition number of the Schur forms for Set2C and Set3C are very similar; this is not surprising because, as we showed in [5] for the CG method, that regardless of which equation set is selected the resulting Schur form pressure equation should represent the same governing dynamics.

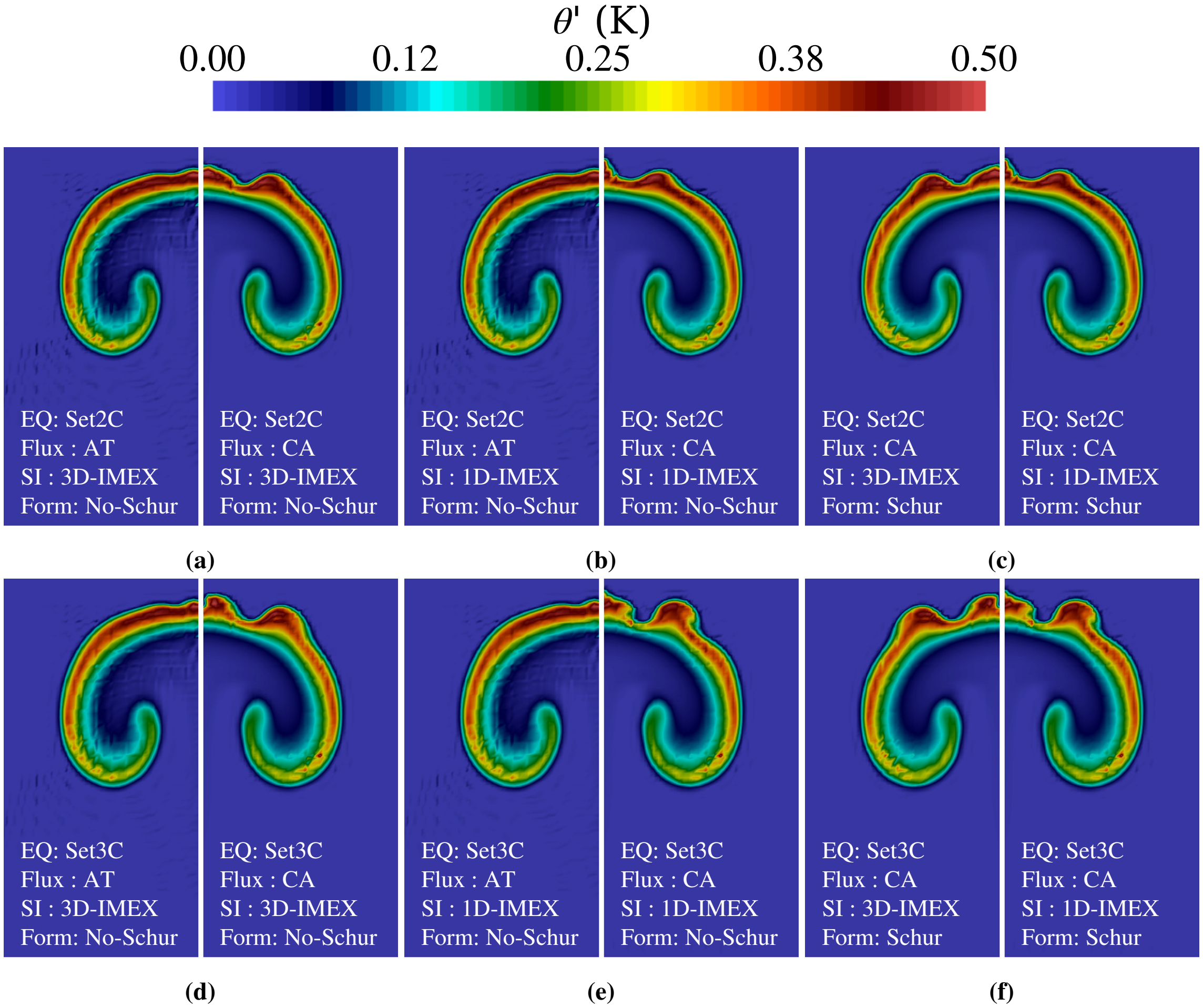

**Figure 2:** Rising Thermal Bubble. The distribution of potential temperature perturbation at $t = 650\ s$ using Set2C with: (a) No-Schur 3D-IMEX, (b) No-Schur 1D-IMEX, and (c) Schur 3D vs. 1D IMEX; Set3C with: (d) No-Schur 3D-IMEX, (e) No-Schur 1D-IMEX, and (f) Schur 3D vs. 1D IMEX.

Table 2 shows the various performance metrics for the rising thermal bubble problem. It shows that the solution metrics (min., max., mass and energy loss) are all in good agreement for various formulations. The IMEX formulations using the *AT* flux are all in agreement with the solution obtained using the fully explicit third-order RK35 (also using the *AT* flux formulation). All equations sets and numerical formulations are able to conserve mass up to machine precision, whereas only Set3C is able to conserve energy up to machine precision. Set2C cannot conserve energy because the thermodynamic equation is written in terms of (density) potential temperature. Similar behaviors were observed when the two equation sets were analyzed using both the continuous Galerkin [23] and discontinuous Galerkin [15] methods. The minimum and maximum values of the thermodynamic variable are consistent across the different IMEX and flux formulations and the mean-absolute-error is low for all cases.

The results in Table 2 show that the 3D-IMEX *No-Schur CA* flux formulation is approximately 1.15 times faster than the *AT* flux formulation for both Set2C and Set3C. This is due to the better conditioning of the system from the *CA* flux formulation, thereby requiring fewer GMRES iteration to reach convergence. The Schur formulation is approximately 1.3 times faster than the *No-Schur AT* flux formulation due to the better conditioning of the system and a smaller overall size of the system, requiring fewer overall FLOPs (floating point operations); this conclusion holds for both Set2C and Set3C. Negligible differences are seen in the computing times between the two *No-Schur* flux formulations for 1D-IMEX. This is because the size of the Jacobian matrices is the same and, therefore, requires the same number of FLOPs for a direct solve (forward reduction and back substitution) of the linear system (standard LU decomposition requires $\mathcal{O}\left(\frac{2}{3}N^3\right)$ operations for $N$ gridpoints, see, e.g., [46]). However, the linear system for the Schur formulation is much smaller (five times smaller) and therefore, the Schur formulation is approximately 4.8 times faster than the 1D-IMEX *No-Schur* discretization with *AT* flux. The Schur form for 1D-IMEX is slightly less than 5 times faster than the *No-Schur* form due to the additional FLOPs that are needed to extract the prognostic variables $(\rho, \mathbf{U}, \Theta, E)$ from the diagnostic variable $(P)$. All timings are for simulations performed using 30 MPI ranks. Note that for this isotropic grid resolution, the 1D-IMEX method is at a disadvantage, which is evident by the large time-to-solution.

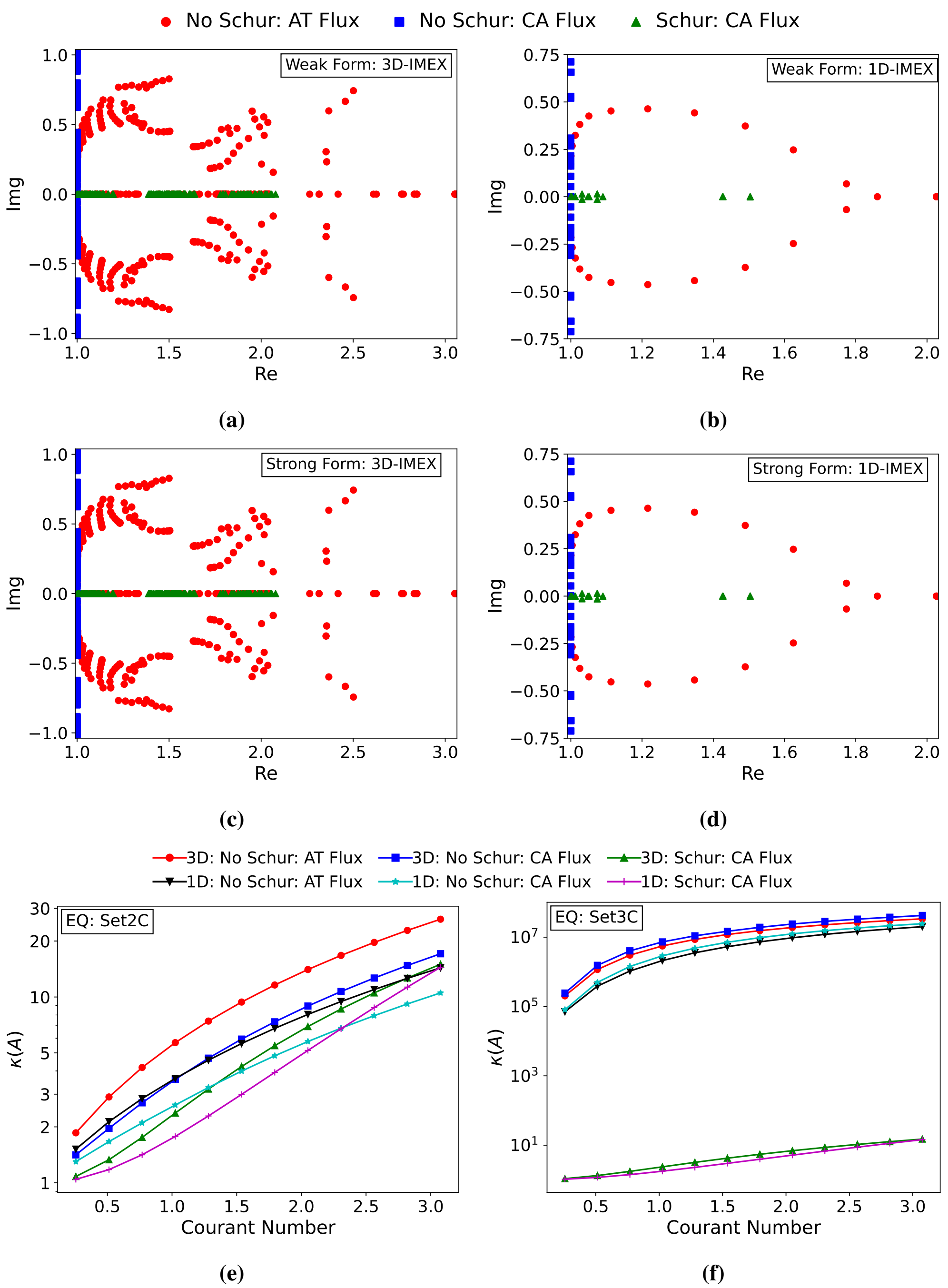

**Figure 3:** Rising Thermal Bubble. The eigenspectra for Set3C for various flux formulations with Courant number $\approx 1$ for variational weak form for: (a) 3D-IMEX and (b) 1D-IMEX; variational strong form for: (c) 3D-IMEX and (d) 1D-IMEX; condition number as a function of Courant number for (e) Set2C and (f) Set3C. (Similar eigenspectra are obtained with Set2C).

**Table 2:** Rising Thermal Bubble. Comparison of performance metrics for various formulations using Set2C and Set3C.

| | | Explicit | 3D-IMEX | | | 1D-IMEX | | |
|---|---|---|---|---|---|---|---|---|
| | | | No-Schur | | Schur | No-Schur | | Schur |
| | | (Rusanov) | *AT* Flux | *CA* Flux | | *AT* Flux | *CA* Flux | |
| | Courant No. | 0.1 | 1.5 | | | 0.1 | | |
| Set2C | Min $\theta'$ | -1.45E-1 | -1.45E-1 | -2.94E-2 | -2.94E-2 | -1.45E-1 | -5.62E-2 | -5.62E-2 |
| | Max $\theta'$ | 5.30E-1 | 5.30E-1 | 5.54E-1 | 5.54E-1 | 5.30E-1 | 5.43E-1 | 5.43E-1 |
| | $\Delta\mathcal{M}$ | 0.0E0 | 0.0E0 | 0.0E0 | 0.0E0 | 2.09E-16 | 0.0E0 | 2.09E-16 |
| | $\Delta\mathcal{E}$ | 2.50E-8 | 2.49E-8 | 2.47E-8 | 2.47E-8 | 2.49E-8 | 2.26E-8 | 2.26E-8 |
| | MAE | - | 1.63E-5 | 2.37E-05 | 2.08E-05 | 9.91E-7 | 2.51E-4 | 2.68E-4 |
| | Sol. Time (s) | 4384 | 319 | 275 | 250 | 15499 | 15580 | 3231 |
| Set3C | Min $\theta'$ | -5.29E-2 | -5.29E-2 | -3.45E-2 | -3.45E-2 | -5.29E-2 | -2.74E-2 | -2.67E-2 |
| | Max $\theta'$ | 5.20E-1 | 5.20E-1 | 5.19E-1 | 5.19E-1 | 5.20E-1 | 5.08E-1 | 5.08E-1 |
| | $\Delta\mathcal{M}$ | 0.0E0 | 0.0E0 | 0.0E0 | 0.0E0 | 2.09E-16 | 2.09E-16 | 2.09E-16 |
| | $\Delta\mathcal{E}$ | 0.0E0 | 0.0E0 | 0.0E0 | 0.0E0 | 0.0E0 | 0.0E0 | 0.0E0 |
| | MAE | - | 1.60E-5 | 2.18E-05 | 1.90E-05 | 6.00E-7 | 3.00E-4 | 3.20E-4 |
| | Sol. Time (s) | 4435 | 322 | 274 | 248 | 15484 | 15443 | 3187 |

## 6.2 Case 2: Density Current

The density current test problem was solved on a $128 \times 32$ mesh of fourth-order elements ($\Delta x = \Delta z = 50\ m$) and is evolved until $t = 900\ s$. As previously mentioned, a uniform diffusion coefficient of $\mu = 75\ m^2 s^{-1}$ is used to maintain stability.

Figure 4 shows the distribution of the potential temperature obtained using the two variational forms and different flux schemes for the 3D-IMEX and 1D-IMEX formulations. Comparing Figs. 4, it can be seen that the different flux formulations have little to no effect on the solution. This is because the dissipation introduced by the artificial viscous operator is significantly greater than the dissipation from the numerical fluxes. The potential temperature obtained using 3D-IMEX, 1D-IMEX in both weak and strong forms are all similar.

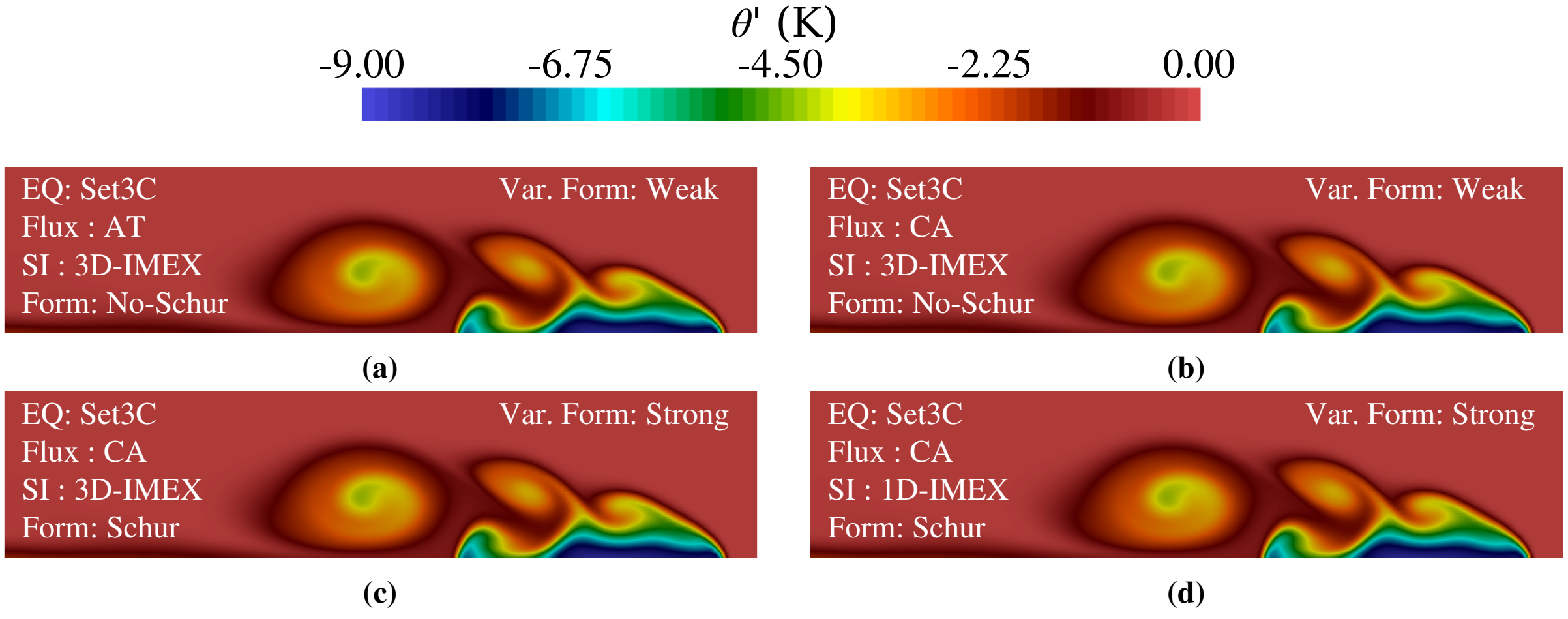

**Figure 4:** Density Current. The distribution of potential temperature perturbation at $t = 900\ s$ using Set3C with variational weak *No-Schur* form and: (a) 3D-IMEX with *AT* flux, (b) 3D-IMEX with *CA* flux; and variational strong Schur form for: (c) 3D IMEX, (d) 1D-IMEX.

Implicit time integration (whether treating the linear or nonlinear terms) introduces greater dissipation into the system than explicit methods [47]. Therefore, the accuracy of the IMEX formulations must be investigated. The accuracy is quantified using the mean absolute error (MAE) of the difference between the solution resulting from the IMEX discretization and some reference solution. The reference solution is computed using a third-order, five-stage, explicit Runge-Kutta method (RK35) [40] using a very small time step (Courant number$\approx 0.0002$) and the *AT* flux formulation, and is compared at $t = 100\ s$. The IMEX schemes are solved using a two-stage, second-order ARK (ARK2) [14], a three-stage, third-order ARK (ARK3) and a four-stage, fourth-order ARK (ARK4) [48] schemes. The different Courant numbers for the 3D-IMEX cases are obtained by modifying the time step appropriately. For the 1D-IMEX cases, the Courant numbers are modified by changing the resolution of the vertical discretization, so that the influence of only vertical/radial operators is affected. Since the grids for the 3D-IMEX and 1D-IMEX cases are different (even for the same Courant number), the errors between the two cannot be compared fairly.

Figure 5 shows the accuracy of the IMEX formulations for the *No-Schur* form of Set3C. The IMEX methods using the *AT* flux formulation exhibit the expected order of convergence relative to the reference solution obtained using RK35 with *AT* flux. For smaller Courant numbers, the IMEX methods using the *CA* flux formulation differ significantly from the reference solution but are in good agreement for larger Courant numbers. This is due to the difference in solution caused by the difference in the amount of dissipation from the numerical flux. For larger Courant numbers, the dissipation from the time integration offsets this difference in dissipation from the numerical fluxes, thereby leading to similar results as the total *AT* flux. In addition, the IMEX methods using the *CA* flux also exhibit the expected order of convergence when the reference solution from the RK35 is obtained using the *CA* flux formulation (Figs. 5c and 5d).

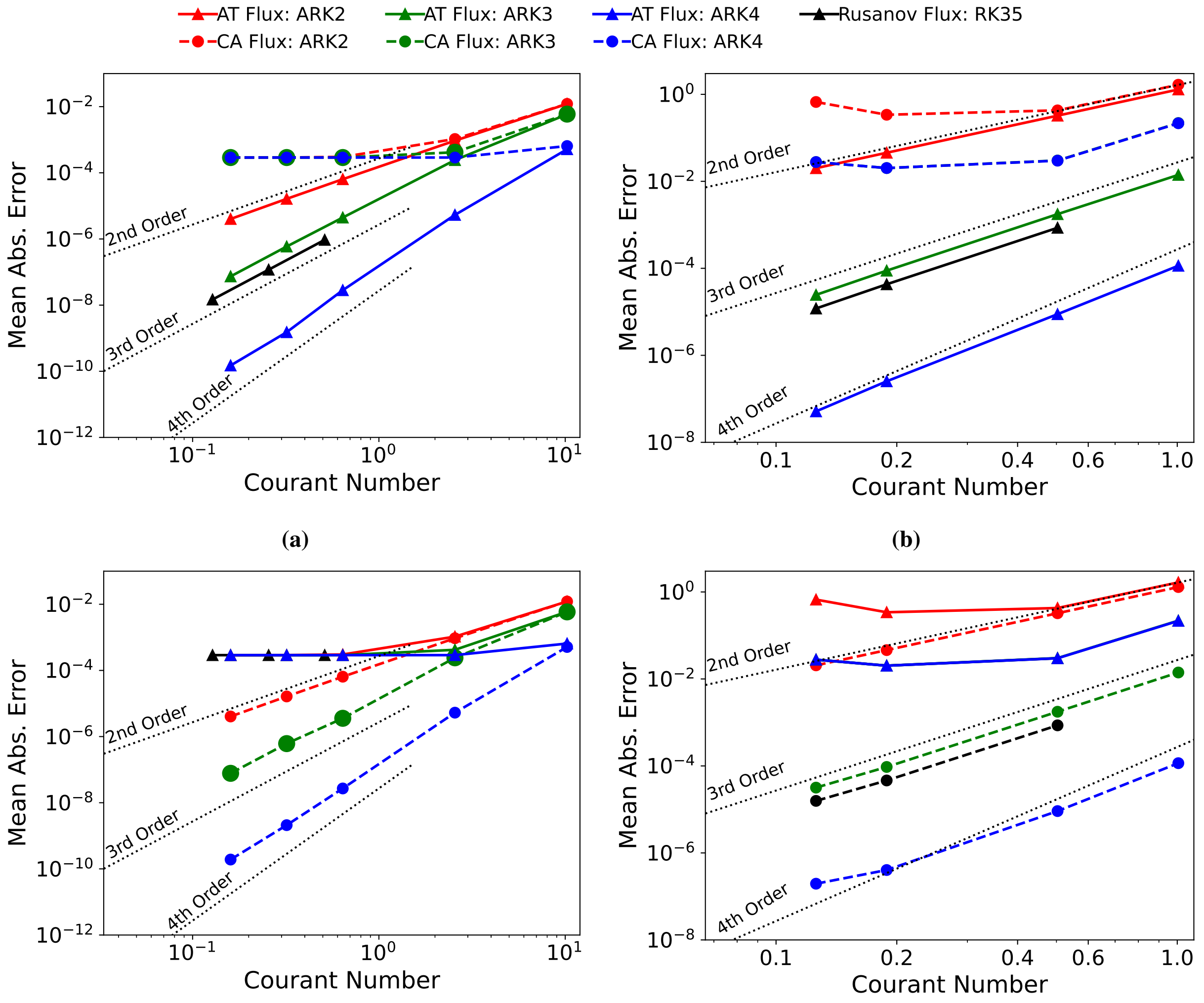

**(a)** **(b)**

**(c)** **(d)**

**Figure 5:** Density Current. The accuracy of the IMEX discretizations and ARK time integrators for Set3C relative to a reference solution $\mathbf{q}^{ref}$ computed using RK35 (with *AT* flux and Courant number $\approx 0.0002$) for: (a) 3D-IMEX and (b) 1D-IMEX formulations; and relative to a reference solution $\mathbf{q}^{ref}$ computed using RK35 (with *CA* flux and Courant number $\approx 0.0002$) for: (c) 3D-IMEX and (d) 1D-IMEX formulations. (Note: the results for 1D-IMEX with *CA* flux is similar for both ARK3 and ARK4).

### 6.3 Case 3: Inertia Gravity Waves

The inertia gravity wave problem was solved on a $96 \times 32$ mesh of fourth-order elements ($\Delta x = 780\ m,\ \Delta z = 78\ m$) and is evolved until $t = 2500\ s$; at this final time, the initial perturbation will be exactly at the midpoint of the domain. The distribution of the potential temperature along the horizontal center line of the domain (at $z = 5km$) obtained using various formulations of the strong form, 3D-IMEX discretization is shown in Fig. 6. No appreciable differences are seen between the distributions from the different discretizations and equation sets. To adhere to the analytical solution, the test case is solved without any artificial diffusion.

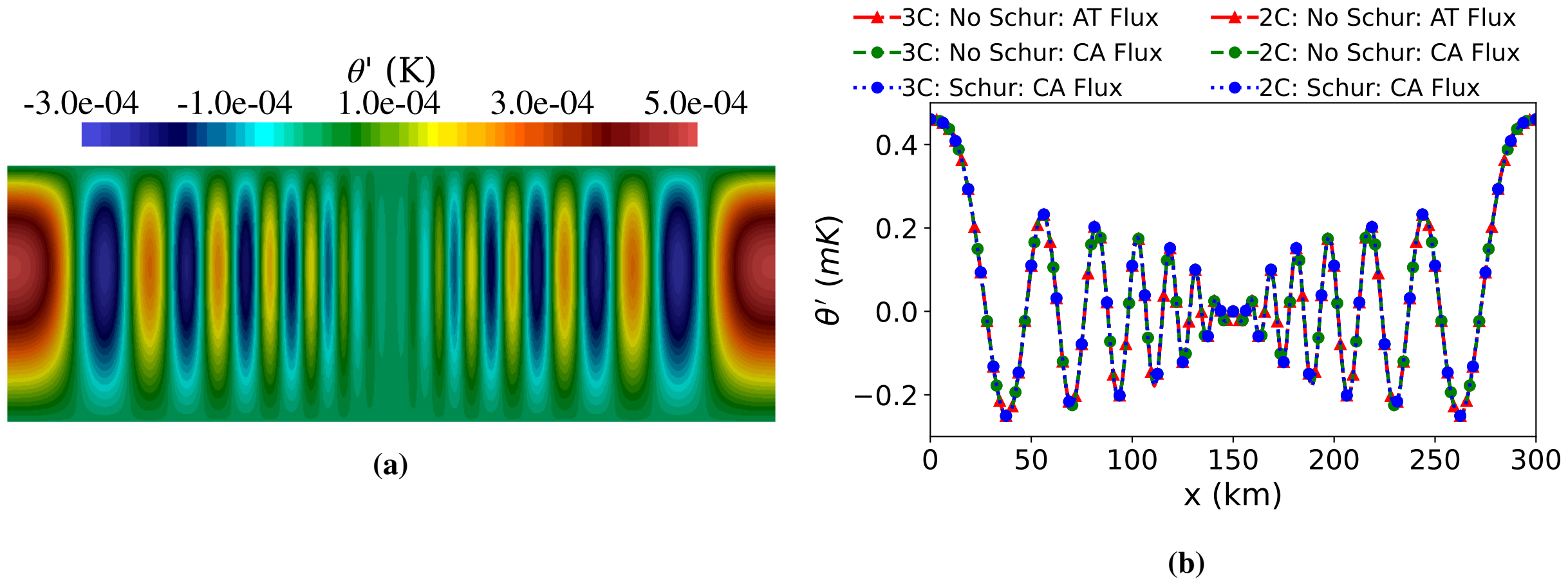

**Figure 6:** Inertia Gravity Wave. The distribution of potential temperature perturbation at $t = 2500\ s$ using 3D-IMEX and: (a) Set2C in No-Schur form with *AT* flux and (b) distribution along the center line for each combination of equation set, numerical flux and IMEX form (note the $y$-axis is in milli-Kelvin).

We now compare the efficiency of the various IMEX formulations. The average number of solver iterations per time step and the time-to-solution (wallclock time) for different Courant numbers are shown in Figs. 7a and 7b, respectively. To not be constrained by the CFL condition on the explicit treatment of the advective dynamics, the constant mean flow is set to $\bar{u} = 0\ m/s$ (see Sec. 5.3). All analyses were performed on 40 cores with ARK2. Therefore, the number of iterations reported per time step are for two ARK stages. The efficiency results presented are shown for both Set2C and Set3C.

For the *No-Schur* system, the *CA* flux formulation requires fewer iterations to reach convergence. This is due to the better conditioning of the system and is consistent with the results for the rising thermal bubble test (see Sec. 6.1). As expected, a larger number of iterations are required with increasing Courant numbers. The Schur system requires fewer iterations than the *No-Schur* system for the same Courant number, due to the better conditioning of the linear system. Although the linear solvers are not the focus of this study, it was previously demonstrated (Fig. 3) that the Schur form results in a symmetric positive definite matrix and can be solved using conjugate gradient (CG). Figure 7a shows that CG requires significantly more iterations than GMRES, particularly for larger Courant numbers. Although CG requires more iterations per time step, the number of FLOPs, storage requirements (no need to store Arnoldi vectors), and inter-processor communication per iteration is much lower than GMRES (when using modified Gram-Schmidt for the Arnolid iteration). Hence, the CG solver is faster for smaller Courant numbers. Previous works have provided conditions for superlinear convergence of the CG [49] and GMRES [50] solvers based on the Ritz values, however, such analysis is beyond the scope of this study but of interest for future work. Both equation sets exhibit similar iteration counts.

Comparing the wallclock time required to reach the simulation time of $500\ s$, Figs. 7b and 7e show that there exists an optimal time step that yields the fastest time-to-solution, after which the efficiency degrades; this behavior can be delayed with preconditioners but we save this for future work (see Sec. 7). For smaller time steps where fewer iterations (one to two) are needed to solve the linear system, the *No-Schur* formulations yield a faster time to solution; similar behavior was seen for continuous Galerkin IMEX formulations [23]. This is due to the computational cost of extracting the prognostic variables from the diagnostic variables being larger than the cost of performing few iterations of the linear solver. At the optimal time step (CN $\approx$ 3), the Schur formulation begins to yield faster time to solution than the *No-Schur* formulations. For Courant numbers smaller than this optimal value, the Schur form with the conjugate gradient solver provides a faster time to solution than the *No-Schur* forms. For smaller Courant numbers (i.e. for more number of time steps), Set3C yields a much faster time to solution than Set2C. This can be attributed to the increased computational cost of the equation of state for Set2C over that of Set3C. For larger Courant numbers (i.e. fewer time steps), the speed-up of Set3C over Set2C is smaller (but Set3C is still faster) due to fewer evaluations of the equation of state.

Figures 7c and 7f show results of the convergence study for the IMEX formulation with the *AT* flux. The reference solution was taken to be the analytical solution of Baldauf and Brdar [37] at $t = 2500\ s$. All methods exhibit the expected order of convergence. A similar breakdown in the order of convergence was seen in [37] and [21]. Figure 7f shows the breakdown of ARK4 at Set3C at higher Courant numbers; this behavior is described in [21, 37].

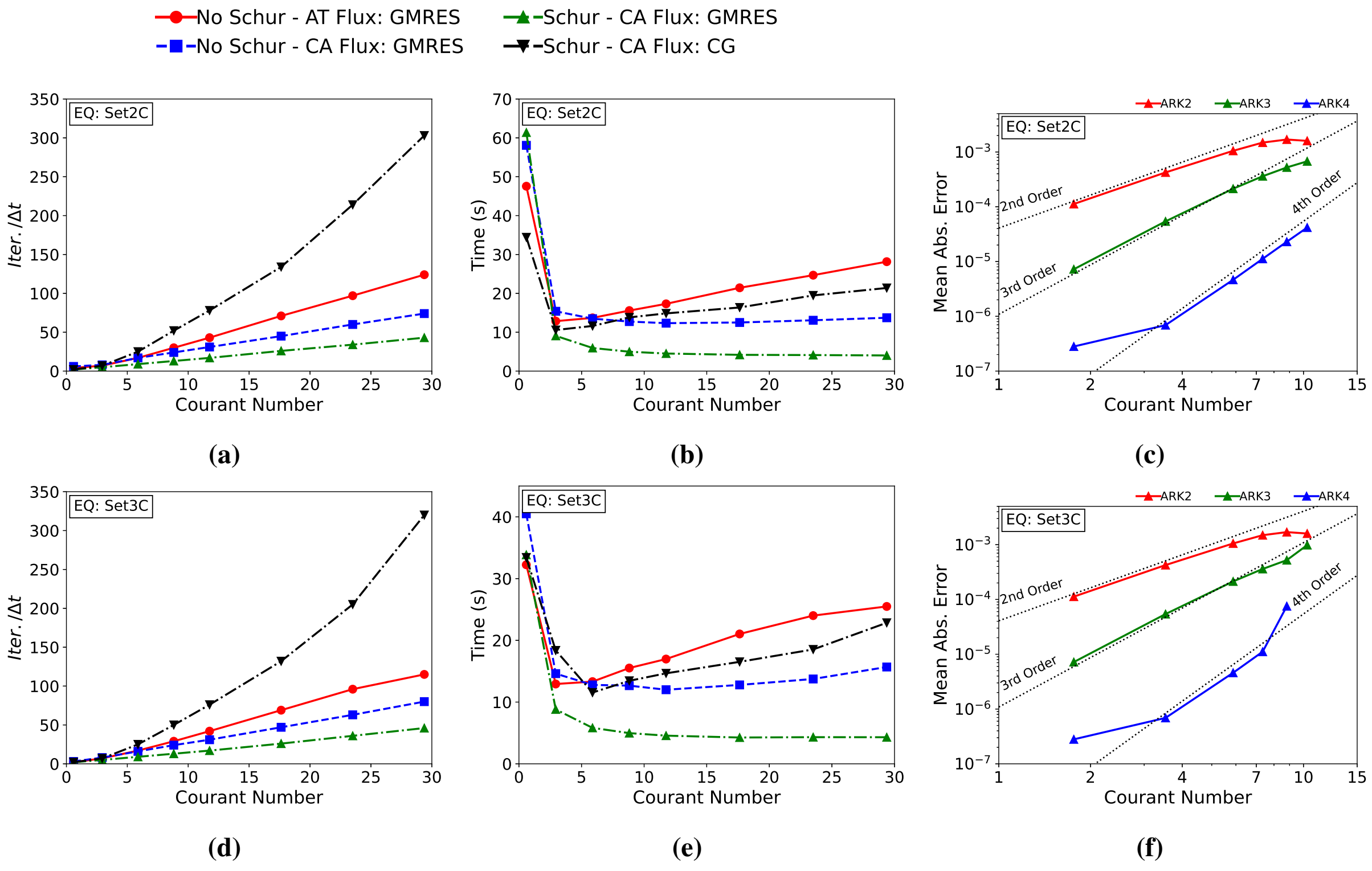

**Figure 7:** Inertia Gravity Wave. Efficiency study of the different equation sets (top row is Set2C and bottom row is Set3C) showing: (a, d) The average number of iterations required to solve the linear system at each time step, (b, e) time required to compute 500 seconds and (c, f) convergence of the various IMEX time integrators with *No-Schur* form and *AT* flux formulation (reference solution was taken to be the analytical solution of Baldauf and Brdar [37]).

### 6.4 Case 4: Acoustic Wave on a Sphere

The acoustic wave on a sphere (AWS) problem is solved on a cubed-sphere grid with a resolution of $\Delta h = 200\ km$ and $\Delta v = 400\ m$, in the horizontal and vertical directions, respectively. The simulation is evolved until $t = 36000s$ with a time step of $\Delta t = 25\ s$ corresponding to Courant number $\approx 37$. The high aspect ratio of such resolutions ($\Delta h/\Delta v = 500$) increases the stiffness due to the vertical resolution and dynamics. This stiffness leads to a drastic increase in the number of iterations required for convergence of the 3D-IMEX system. Therefore, the 1D-IMEX discretization with direct solver is ideal for such applications. Figure 8 shows snapshots at different times for the acoustic wave traveling around the sphere. The solutions are obtained using the 1D-IMEX discretization of Set3C with the *AT* flux formulation and ARK2 time integration. No significant differences are seen between the *AT* and *CA* flux formulations.

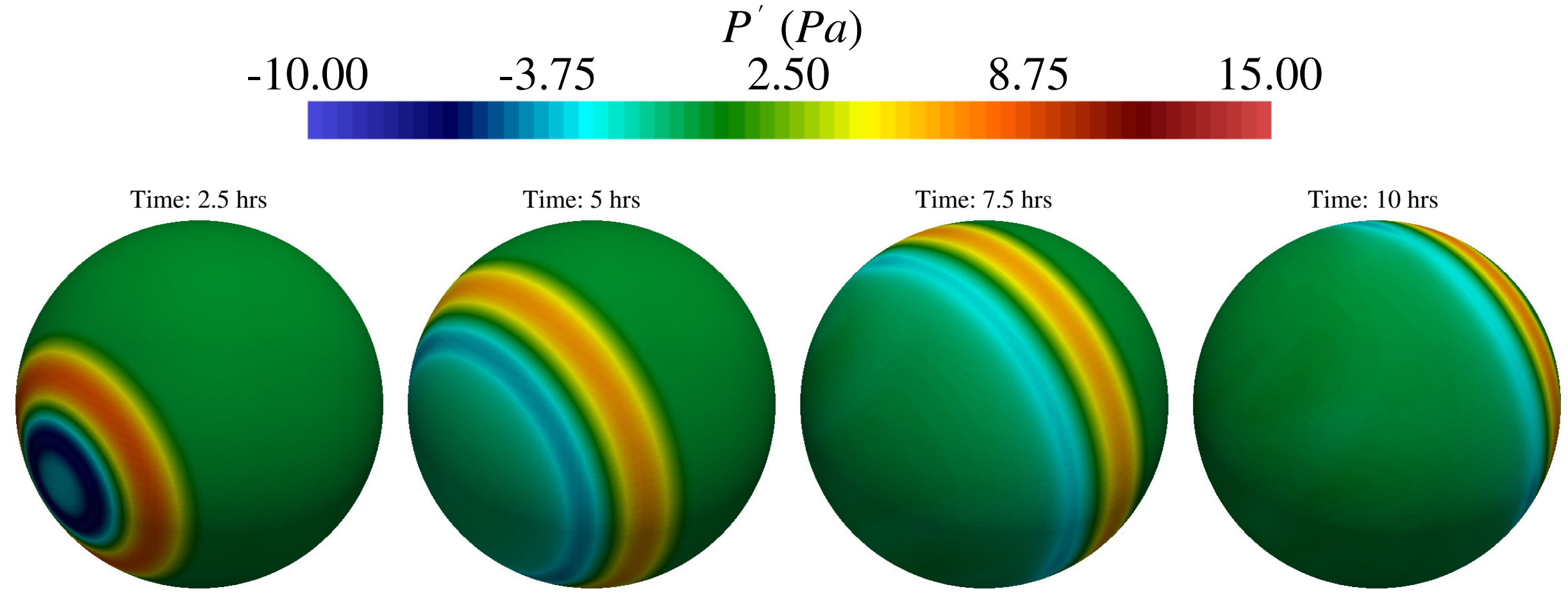

**Figure 8:** Acoustic Wave on Sphere. Pressure perturbation at different times obtained using 1D-IMEX discretization of Set3C with the No-Schur *AT* flux formulation.

The efficiency of the 1D-IMEX over the 3D-IMEX for global problems (i.e., problems on the sphere) is investigated. The number of iterations and time-to-solution for different Courant numbers is investigated for all formulations presented. Since the linear system in the 1D-IMEX formulation is solved using a direct method (banded LU decomposition from the LAPACK library), studying its performance by simply changing the time steps to obtain the required Courant number is not appropriate (since the matrix size and FLOPs will remain the same and, therefore, the performance will increase linearly with the increase in time step)).

For this reason, the Courant number in the efficiency study of the 1D-IMEX formulations is updated by modifying the vertical resolution. That is, higher resolutions yield larger Courant numbers and larger system matrices. For the 3D-IMEX case the grid resolution is kept constant; instead the required Courant number is attained by adapting the time step. This keeps the size of the linear system constant between different Courant numbers but with different stiffness, thereby requiring different number of iterations to reach convergence. All analyses are performed with Set3C and the ARK2 scheme using 10 MPI ranks.

Figure 9 shows the results of the convergence study for 1D-IMEX (with *AT* flux), average number of iterations (for the 3D-IMEX case) and the time required to evolve the system for 10 time steps. The Courant number presented is the Courant number of the acoustic term with vertical grid spacing (i.e., $a = \sqrt{\frac{\gamma P}{\rho}}$ and $\Delta s = \Delta v$ in Eq. (5.4)). Comparing Figs. 9a and 9b shows that the 1D-IMEX methods are more than an order of magnitude faster than the 3D-IMEX methods for such high-aspect ratio problems. This is because the 3D-IMEX methods are restricted by the stiffness due to the vertical discretization and dynamics. The Schur forms for both 3D-IMEX and 1D-IMEX are significantly faster than their *No-Schur* forms. In the case of 1D-IMEX, the timings of the two *No-Schur* forms coincide as expected, since the two linear systems are of the same size and, because we are using direct solvers, the condition number of the matrices play no role with respect to time-to-solution. Although the 3D-IMEX formulations allow for a larger *stable* time step (since acoustic terms in all directions are handled implicitly), the larger computing/simulation times often render these methods unusable for such global models. Hence, the 1D-IMEX are ideally suited for high-aspect ratio meshes as they allow for much larger *usable* time steps while greatly reducing the time to solution.

Figure 9c shows the average number of iterations per time step required to reach convergence. The Schur form of the 3D-IMEX discretization requires significantly fewer iterations to reach convergence; this can be attributed to the better conditioning (smaller condition number) as well as the location of the eigenvalues of the Schur system (Krylov methods prefer the eigenvalues near the real axis). The *No-Schur* form with *AT* and *CA* fluxes require similar number of iterations for smaller Courant numbers (CN $< 5$). For larger Courant numbers the *CA* flux formulation requires fewer iterations than the *AT* flux, thereby yielding a much faster time to solution. Due to the computational cost of the 3D-IMEX formulations for such highly anisotropic grids, the convergence analysis is only performed for 1D-IMEX formulations. These 1D-IMEX schemes (Fig. 9d) are shown to follow the expected order of convergence.

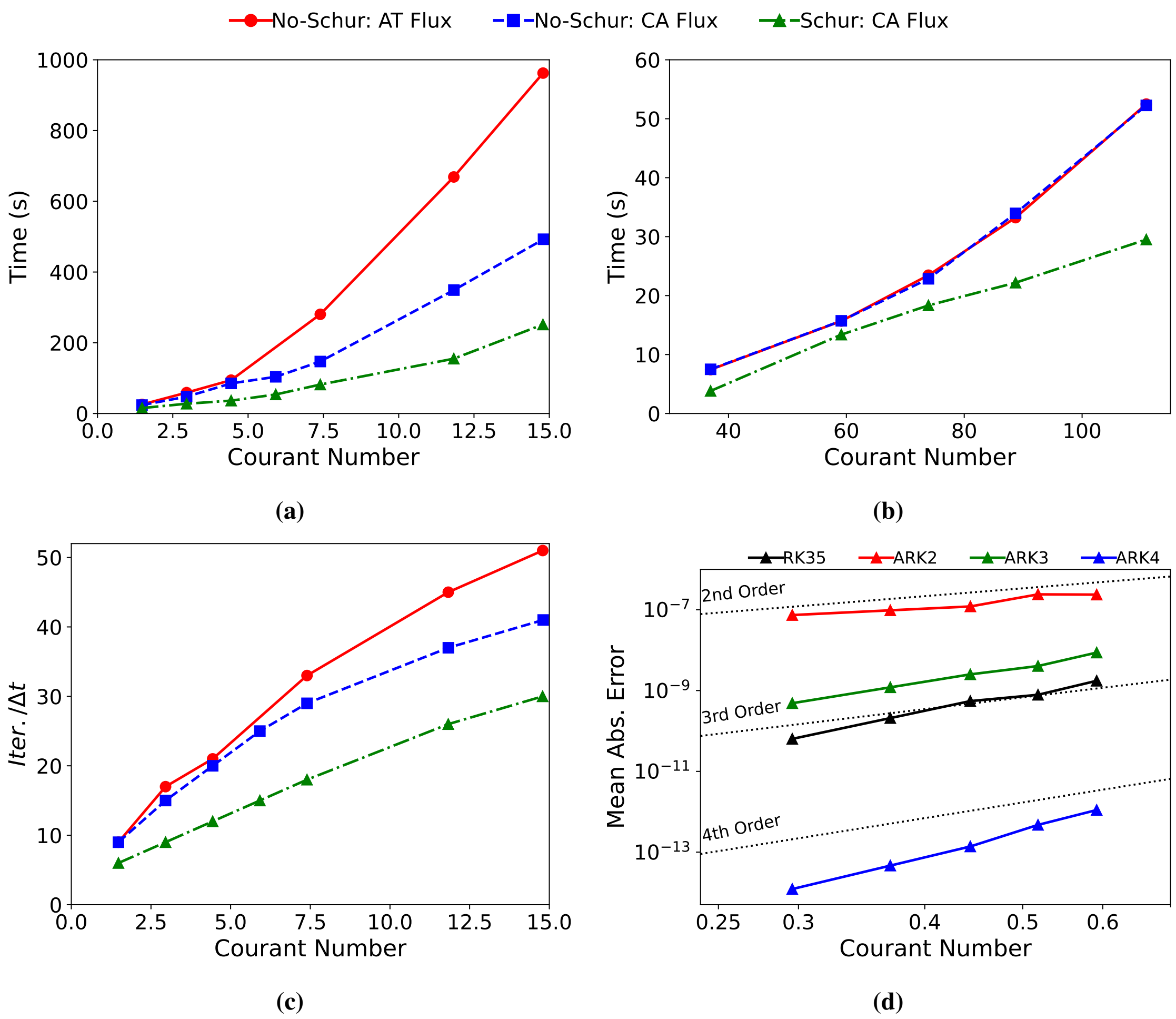

**Figure 9:** Acoustic Wave on Sphere. The accuracy and efficiency of the IMEX discretizations with Set3C showing: the time required to compute 10 time steps as a function of Courant number using: (a) 3D-IMEX and (b) 1D-IMEX; (c) average number of iterations per time step for 3D-IMEX and (d) convergence of the 1D-IMEX scheme (reference solution was obtained using RK35 and CN $\approx 0.001$)

# 7 Conclusion

We present and analyze IMplicit-EXplicit formulations for discontinuous Galerkin discretizations of different forms of the Euler equations. Two different IMEX methods that resolve the stiffness due to the governing dynamics and geometric discretization are shown.The 3D-IMEX method, which treats the acoustical stiffness in all directions, allows for significantly larger time steps and yields a faster time to solution for problems with mesh aspect ratios closer to one. In the case of global models, where the aspect ratio is much larger, the 1D-IMEX method is more efficient and more than an order of magnitude faster than 3D-IMEX. For each IMEX formulation and governing equation set, efficient Schur complements are derived to reduce the size of the linear system from a system of five variables to a single equation for pressure. The derivation and consequences of the different numerical fluxes used in these Schur complements are discussed. The Schur formulations of these IMEX methods result in a smaller, better conditioned system, and, therefore, require few iterations to reach convergence. Moreover, the resulting linear system is shown to have purely real eigenvalues confirming the ellipticity of the resulting discrete form of the Helmholtz-like differential operator; this last point is important because it now opens the door to developing preconditioners that will further improve the performance of this method (e.g., multigrid methods). The ARK methods of different orders also yield the expected theoretical rates of convergence for the test cases considered. These different formulations of the IMEX methods are shown to be efficient and accurate for both mesoscale (box) and global (sphere) models.

## Acknowledgments

The authors gratefully acknowledge the support by the generosity of Eric and Wendy Schmidt by recommendation of the Schmidt Futures program. We are also grateful for funding provided by the Office of Naval Research under grant # N0001419WX00721 and National Science Foundation under grant AGS-1835881. This work was performed in part when S. Reddy, M. Waruszewski, and F. A.V. B. Alves held National Academy of Sciences' National Research Council (NRC) Fellowships at the Naval Postgraduate School. FXG would like to thank Prof. Maria Lukacova-Medvidova for our conversation on this topic at the Moist Processes in the Atmosphere meeting in Oberwolfach which formed the initial impetus for this work.

## References

[1] Francis X. Giraldo. *An Introduction to Element-Based Galerkin Methods on Tensor-Product Bases - Analysis, Algorithms and Applications*. Springer, 2020.

[2] L. Yelash, A. Müller, M. Lukacova-Medvid'ova, F. X. Giraldo, and V. Wirth. Adaptive discontinuous evolution Galerkin method for dry atmospheric flow. *Journal of Computational Physics*, 268:106–133, JUL 1 2014.

[3] V. Dolejší and M. Feistauer. A semi-implicit discontinuous Galerkin finite element method for the numerical solution of inviscid compressible flow. *Journal of Computational Physics*, 198(2):727–746, Aug. 2004.

[4] Alex Kanevsky, Mark H. Carpenter, David Gottlieb, and Jan S. Hesthaven. Application of implicit-explicit high order runge-kutta methods to discontinuous-galerkin schemes. *Journal of Computational Physics*, 225(2):1753–1781, Aug. 2007.

[5] F. X. Giraldo and M. Restelli. High-order semi-implicit time-integrators for a triangular discontinuous Galerkin oceanic shallow water model. *International Journal for Numerical Methods in Fluids*, 63(9):1077–1102, Jul. 2010.

[6] M. Restelli and F. X. Giraldo. A conservative discontinuous galerkin semi-implicit formulation for the navier-stokes equations in nonhydrostatic mesoscale modeling. *SIAM Journal on Scientific Computing*, 2009.

[7] M. Kwizak and A.J. Robert. A semi–implicit scheme for grid point atmospheric models of the primitive equations. *Monthly Weather Review*, 99(1):32–36, January 1971.

[8] A.J. Simmons, B.J. Hoskins, and D.M. Burridge. Stability of the semi-implicit method of time integration. *Monthly Weather Review*, 106(3):405–412, March 1978.

[9] Klaus Kaiser, Jochen Schutz, Ruth Schoebel, and Sebastian Noelle. A new stable splitting for the isentropic euler equations. *Journal of Scientific Computing*, 70(3):1390–1407, Mar 2017.

[10] Klaus Kaiser and Jochen Schutz. A high-order method for weakly compressible flows. *Communications in Computational Physics*, 22(4):1150–1174, Oct. 2017.

[11] Georgij Bispen, Maria Lukacova-Medvid'ova, and Leonid Yelash. Asymptotic preserving imex finite volume schemes for low mach number euler equations with gravitation. *Journal of Computational Physics*, 335:222–248, Apr. 2017.

[12] Jonas Zeifang, Klaus Kaiser, Andrea Beck, Jochen Schuetz, and Claus-Dieter Munz. Efficient high-order discontinuous galerkin computations of low mach number flows. *Communications in Applied Mathematics and Computational Science*, 13(2):243–270, 2018.

[13] Jonas Zeifang, Jochen Schuetz, Klaus Kaiser, Andrea Beck, Maria Lukacova-Medvid'ova, and Sebastian Noelle. A novel full-euler low mach mumber imex splitting. *Communications in Computational Physics*, 27(1):292–320, Jan 2020.

[14] Francis X. Giraldo, James F. Kelly, and Emil M. Constantinescu. Implicit-explicit formulations of a three-dimensional nonhydrostatic unified model of the atmosphere (NUMA). *SIAM Journal of Scientific Computing*, 35(5):B1162–B1194, 2013.

[15] Francis X. Giraldo and Marco Restelli. A study of spectral element and discontinuous galerkin methods for the navier-stokes equations in nonhydrostatic mesoscale atmospheric modeling: Equation sets and test cases. *Journal of Computational Physics*, 227:3849–3877, 2008.

[16] William C. Skamarock, Joseph B. Klemp, Jimy Dudhia, David O. Gill, Zhiquan Liu, Judith Berner, Wei Wang, Jordan G. Powers, Michael G. Duda, Dale Barker, and Xiang yu Huang. A description of the advanced resaerch WRF model version 4.3. Technical report, NCAR, 2021.

[17] Nash'at Ahmad and John Lindeman. Euler solutions using flux-based wave decomposition. *International Journal for Numerical Methods in Fluids*, 54:41–72, 2007.

[18] Robert V. Rohli and Chunyan Li. *Meteorology for Coastal Scientists*. Springer, 2021.

[19] John. D. Anderson. Governing equations of fluid dynamics. In John F. Wendt, editor, *Computational Fluid Dynamics*, chapter 2, pages 15–51. Springer, Berlin, Heidelberg, 1992.

[20] Akshay Sridhar, Yassine Tissaoui, Simone Marras, Zhaoyi Shen, Simone Byrne Charles Kawczynski, Kiran Pamnany, Maciej Waruzewski, Jeremy E. Kozdon Thomas H. Gibson, Valentin Churavy, Lucas C. Wilcox, Francis X. Giraldo, and Tapio Schneider. Large-eddy simulations with ClimateMachine v0.2.0: A new open-source code for atmospheric simulations on GPUs and CPUs. *Geoscientific Model Development*, Preprint:1–41, 2021.

[21] Maciej Waruszewski, Jeremy E. Kozdon, Lucas C. Wilcox, Thomas H. Gibson, and Francis X. Giraldo. Entropy stable discontinuous galerkin methods for balance laws in non-conservative form: Applications to euler with gravity. *Journal of Computational Physics*, 468:111507, 2022.

[22] Gregor J. Gassner, Andrew R. Winters, and David A. Kopriva. Split form nodal discontinuous galerkin schemes with summation-by-parts property for the compressible euler equations. *Journal of Computational Physics*, 327:39–66, 2016.

[23] Francis X. Giraldo, Marco Restelli, and M. Lauter. Semi-implicit formuations of the navier-stokes equations: Application ot nonhydrostatic atmospheric modeling. *SIAM Journal of Scientific Computing*, 32(6):3394–3425, 2010.

[24] Bernardo Cockburn and Chi-Wang Shu. The local discontinuous galerkin method for time-dependent convective-diffusion systems. *SIAM Journal of Numerical Analysis*, 35(6):2440–2463, 1998.

[25] Sohail Reddy, Yassine Tissaoui, Felipe A. V. de Branganca Alves, Simone Marras, and Francis X. Giraldo. Comparison of sub-grid scale models for large-eddy simulation using a high-order spectral element approximation of the compressible navier-stokes equations at low mach number. *arXiv*, 2021.

[26] S. Marras, M. Nazarov, and F. X. Giraldo. Stabilized high-order galerkin methods based on a parameter-free dynamic SGS model for LES. *Journal of Computational Physics*, 301:77–101, 2015.

[27] A. W. Vreman. An eddy-viscosity subgrid-scale model for turbulent shear flow: Algebraic theory and applications. *Physics of Fluids*, 16(10):3670–3681, 2004.

[28] David L. Williamson. The evolution of dynamical cores for global atmospheric models. *The Journal of the Meteorological Society of Japan*, 85B:241–269, 2007.

[29] Andreas Muller, Michal A. Kopera, Simone Marras, Lucas C. Wilcox, Tobin Isaac, and Francis X. Giraldo. Strong scaling for numerical weather prediction at petascale with the atmospheric model NUMA. *The International Journal of High Performance Computing Applications*, 33(2), 2019.

[30] Daniel S. Abdi and Francis X. Giraldo. Efficient construction of unified continuous and discontinuous galerkin formulations for the 3d euler equations. *Journal of Computational Physics*, 320:46–68, 2016.

[31] Daniel S. Abdi, Lucas C. Wilcox, Timothy C. Warburton, and Francis X. Giraldo. A GPU-accelerated continuous and discontinuous Galerkin non-hydrostatic atmospheric model. *International Journal of High Performance Computing Applications*, 33(1):81–109, 2017.

[32] Daniel S. Abdi, Francis X. Giraldo, Emil M. Constantinescu, Lester E. Carr, Lucas C. Wilcox, and Timothy C. Warburton. Acceleration of the implicit-explicit non-hydrostatic unified model of the atmosphere (NUMA) on manycore processors. *The International Journal of High Performance Computing Applications*, 33(2):242–267, 2017.

[33] Carsten Burstedde, Lucas C. Wilcox, and Omar Ghattas. `p4est`: Scalable algorithms for parallel adaptive mesh refinement on forests of octrees. *SIAM Journal on Scientific Computing*, 33(3):1103–1133, 2011.

[34] Rao Kotamarthi, Katharine Hayhoe, Linda O. Mearns, Donald Wuebbles, Jennifer Jacobs, and Jennifer Jurado. *Downscaling Techniques for High-Resolution Climate Projections: From Global Change to Local Impacts*. Cambridge University Press, Cambridge, 2021.

[35] Andre Robert. Bubble convection experiments with a semi-implicit formulation of the euler equations. *Journal of Atmospheric Sciences*, 50(13):1865–1873, 1993.

[36] J. M. Straka, Robert B. Wilhelmson, Louis J. Wicker, John R. Anderson, and Kevin K. Droegemeier. Numerical solutions of a non-linear density current: A benchmark solution and comparisons. *International Journal for Numerical Methods in Fluids*, 17(1):1–22, 1993.

[37] Michael Baldauf and Slavko Brdar. An analytic solution for linear gravity waves in a channel as a test for numerical models using the non-hydrostatic, compressible euler equations. *Quarterly Journal of the Royal Meteorological Society*, 139(677):1977–1989, 2013.

[38] William C. Skamarock and Joseph B. Klemp. Efficiency and accuracy of the klemp-wilhelmson time-splitting technique. *Monthly Weather Review*, 122(11):2623–2630, 1994.

[39] Hirofumi Tomita and Masaki Satoh. A new dynamical framework of nonhydrostatic global model using the icosahedral grid. *Fluid Dynamics Research*, 34(6):357–400, 2004.

[40] Raymond J. Spiteri and Steven J. Ruuth. A new class of optimal high-order strong-stability-preserving time discretization methods. *SIAM Journal on Numerical Analysis*, 40(2), 2002.

[41] Thomas M. Bendall, Thomas H. Gibson, Jemma Shipton, and Colin J. Cotter. A compatible finite-element discretization for the moist compressible euler equations. *Quarterly Journal of the Royal Meteorological Society*, 146:3187–3205, 2019.

[42] L.E. Carr, C.F. Borges, and Francis X. Giraldo. Matrix-free polynomial-based nonlinear least squares optimized preconditioning and its application to discontinuous galerkin discretizations of the euler equations. *Jounral of Scientific Computing*, 66:917–940, 2016.

[43] Robert M. Kirby and George Em Karniadakis. Selecting the numerical flux in discontinuous galerkin methods for diffusion problems. *Journal of Scientific Computing*, 22 - 23, 2005.

[44] F. Bassi and S. Rebay. A high-order accurate discontinuous finite element method for the numerical solution of the compressible navier-stokes equations. *Journal of Computational Physics*, 131(2):267–279, 1997.

[45] Carlos Erik Baumann and J. Tinsley Oden. A discontinuous *hp* finite element method for convection-diffusion problems. *Computer Methods in Applied Mechanics and Engineering*, 175(3-4):311–341, 1999.

[46] L.N. Trefethen and D. Bau III. *Numerical Linear Algebra*. SIAM, Philadelphia, 1997.

[47] Philipp Öffner, Jan Glaubitz, and Hendrik Ranocha. Analysis of artificial dissipation of explicit and implicit time-integration methods. *arXiv*, 2016.

[48] Christopher A. Kennedy and Mark H. Carpenter. Additive runge-kutta schemes for convection-diffusion-reaction equations. *Applied Numerical Mathematics*, 44(1-2), 2003.

[49] A. van der Sluis and H.A. van der Vorst. The rate of convergence of conjugate gradient. *Numerische Mathematik*, 48:543–560, 1986.

[50] H.A. van der Vorst and C. Vuik. The superlinear convergence behaviour of gmres. *Journal of Computational and Applied Mathematics*, 48:327–341, 1993.

[51] G. Tumolo, L. Bonaventura, and M. Restelli. A semi-implicit, semi-Lagrangian, p-adaptive discontinuous Galerkin method for the shallow water equations. *Journal of Computational Physics*, 232(1):46–67, JAN 2013.

[52] Giovanni Tumolo and Luca Bonaventura. A semi-implicit, semi-Lagrangian discontinuous Galerkin framework for adaptive numerical weather prediction. *Quarterly Journal of the Royal Meteorological Society*, 141(692, A):2582–2601, OCT 2015.

[53] F. X. Giraldo, J. S. Hesthaven, and T. Warburton. Nodal high-order discontinuous Galerkin methods for the spherical shallow water equations. *Journal of Computational Physics*, 181(2):499–525, September 2002.

# A Effects of Wave Speeds on Flux Formulations for IMEX Discretizations

Here, we present conditions on the wave speed of the total flux, given the treatment of the linear flux so that we obtain the appropriate eigenvalue/wave speed for the nonlinear flux. Computing the nonlinear flux in Eq. (4.19) as $\mathbf{F}_{NL} = \mathbf{F}_T - \mathbf{F}_L$ and writing $\lambda_{NL} = \lambda_T - \lambda_L$, let us consider the following cases:

**Case 1** Treating the implicit system with a centered flux ($\lambda_L = 0$), letting $\lambda_{NL} = u$, the nonlinear (explicit) operator in Eq. (4.19) can be written as

$$\begin{aligned}\int_\Omega \psi \nabla \cdot \mathbf{F}_{NL}\, d\Omega = &\int_\Gamma \psi \left(\{\{\mathbf{F}_T\}\} - \frac{u}{2}[\![\mathbf{q}]\!]\hat{\mathbf{n}}\right) \cdot \hat{\mathbf{n}}\, d\Gamma - \int_\Omega \nabla\psi \cdot \mathbf{F}_T\, d\Omega \\ &- \int_\Gamma \psi \{\{\mathbf{F}_L\}\} \cdot \hat{\mathbf{n}}\, d\Gamma + \int_\Omega \nabla\psi \cdot \mathbf{F}_L\, d\Omega\end{aligned} \tag{A.1}$$

where we obtain a form that treats the total flux (and the nonlinear flux) using Rusanov with the correct (advective) wave speed of the nonlinear flux Jacobian. This shows that when the implicit term is treated with centered fluxes, the total flux $\mathbf{F}_T$ must be treated with a Rusanov flux with wave speed $u$. This analysis shows why the splitting used in [51, 52] works in a semi-Lagrangian setting since the nonlinear (explicit) operator is treated in a semi-Lagrangian approach whereby only the advective speed is used in the method of characteristics and the remainder of the terms used centered fluxes, which is what we have analyzed here. Regardless of how we view this splitting, it is clear that this approach results in a complete system that is less dissipative when $u \ll a$.

**Case 2** Treating the total flux with a linearized Rusanov flux ($\lambda_T = u + a_0$) with $a_0 = \sqrt{\frac{\gamma P_0}{\rho_0}}$, we get

$$\begin{aligned}\int_\Omega \psi \nabla \cdot \mathbf{F}_{NL}\, d\Omega = &\int_\Gamma \psi \left( \{\{\mathbf{F}_T\}\} - \frac{u + a_0}{2} [\![\mathbf{q}]\!] \hat{\mathbf{n}} \right) \cdot \hat{\mathbf{n}}\, d\Gamma - \int_\Omega \nabla\psi \cdot \mathbf{F}_T\, d\Omega \\ &- \int_\Gamma \psi \left( \{\{\mathbf{F}_L\}\} - \frac{\lambda_L}{2} [\![\mathbf{q}]\!] \hat{\mathbf{n}} \right) \cdot \hat{\mathbf{n}}\, d\Gamma + \int_\Omega \nabla\psi \cdot \mathbf{F}_L\, d\Omega \end{aligned} \tag{A.2}$$

where $\lambda_{NL} = u$ only when $\lambda_L = a_0$. This shows that when the implicit term is treated with a linearized Rusanov flux with wave speed $\lambda_L$, the total flux $\mathbf{F}_T$ should be treated with a linearized Rusanov flux with $\lambda_T = u + \lambda_L$. When $\lambda_L = 0$ (as it is in the Schur form), then $\lambda_{NL} = u + a_0$ is the wave speed of the terms being treated explicitly. Hence, the explicit time integration is still subject to the CFL restrictions due to the acoustic dynamics which requires using a smaller time steps.

**Case 3** Treating the total flux with the total nonlinear Rusanov flux ($\lambda_T = u + a$) we get

$$\begin{aligned}\int_\Omega \psi \nabla \cdot \mathbf{F}_{NL}\, d\Omega = &\int_\Gamma \psi \left( \{\{\mathbf{F}\}\} - \frac{u + a}{2} [\![\mathbf{q}]\!] \hat{\mathbf{n}} \right) \cdot \hat{\mathbf{n}}\, d\Gamma - \int_\Omega \nabla\psi \cdot \mathbf{F}\, d\Omega \\ &- \int_\Gamma \psi \left( \{\{\mathbf{F}_L\}\} - \frac{\lambda_L}{2} [\![\mathbf{q}]\!] \hat{\mathbf{n}} \right) \cdot \hat{\mathbf{n}}\, d\Gamma + \int_\Omega \nabla\psi \cdot \mathbf{F}_L\, d\Omega. \end{aligned} \tag{A.3}$$

This appears valid but the pressure $P$ appears nonlinearly in $a$. When $\lambda_L = 0$ (i.e., implicit term is treated with a centered flux), then the wave speed of the nonlinear term is $\lambda_{NL} = u + a \gg u$ which is not an eigenvalue of the nonlinear flux Jacobian and can lead to a large portion of the eigenspectra to lie outside the region of stability of the explicit time integrator thereby requiring a smaller time step that is dictated by the acoustic dynamics. When $\lambda_L = a_0$ (i.e., implicit term is treated with a linearized Rusanov flux), then $\lambda_{NL} = u + a - a_0 = u + a'$ is still not an eigenvalue of the nonlinear operator (as shown above) but yields eigenvalues nearer to it. Linearizing $a$ about $P_0$ and $\rho_0$, we obtain an expression for $a'$.

$$\begin{aligned} a &\approx \underbrace{\sqrt{\frac{\gamma P_0}{\rho_0}}}_{a_0} + \underbrace{\frac{1}{2}\left( \sqrt{\frac{\gamma}{P_0 \rho_0}} P' - \sqrt{\frac{\gamma P_0}{\rho_0^3}} \rho' \right)}_{a'} \\ &\approx a_0 + \frac{a_0}{2} \left( \frac{P'}{P_0} - \frac{\rho'}{\rho_0} \right) \end{aligned}$$

where $P' \ll P_0$ and $\rho' \ll \rho_0$. We can then perform a scale analysis with $\mathcal{O}(\gamma) = \mathcal{O}(\rho_0) = \mathcal{O}(1)$, $\mathcal{O}(P_0) = \mathcal{O}(10^5)$. Figure 8 shows that $\mathcal{O}(P') = \mathcal{O}(10)$, and, numerically we can show that, $\mathcal{O}(\rho') < \mathcal{O}(10^{-4})$. Therefore $\mathcal{O}(a') < \mathcal{O}(10^{-1}) < \mathcal{O}(u)$, hence such a choice of the numerical flux is not problematic for low Mach number applications.

# B The Necessity for Special Choices on the Flux to Construct the Schur Complement

To understand the challenge of constructing the Schur complement for the DG method, let us consider the linearized shallow water equations

$$\begin{aligned} \frac{\partial \varphi}{\partial t} + \nabla \cdot \mathbf{U} &= 0 \\ \frac{\partial \mathbf{U}}{\partial t} + \nabla \varphi &= 0 \end{aligned} \tag{B.1}$$

where $(\varphi, \mathbf{U})$ represents the geopotential height and momentum, respectively. Applying the element-based Galerkin method to Eq. (B.1) (multiplying by a test function and integrating) yields

$$\begin{aligned}\int_{\Omega_e} \psi \frac{\partial \varphi_N}{\partial t}\, d\Omega_e + \int_{\Omega_e} \psi \nabla \cdot \mathbf{U}_N \, d\Omega_e = 0 \\ \int_{\Omega_e} \psi \frac{\partial \mathbf{U}_N}{\partial t}\, d\Omega_e + \int_{\Omega_e} \psi \nabla \varphi_N \, d\Omega_e = 0.\end{aligned} \tag{B.2}$$

Next, note that the second terms for each equation can be written as follows

$$\begin{aligned}\int_{\Omega_e} \psi \nabla \cdot \mathbf{U}_N \, d\Omega_e = \int_{\Gamma_e} \psi \widetilde{\mathbf{U}}_N \cdot \hat{\mathbf{n}}\, d\Gamma_e - \int_{\Omega_e} \nabla \psi \cdot \mathbf{U}_N^{(e)}\, d\Omega_e \\ \int_{\Omega_e} \psi \nabla \varphi_N \, d\Omega_e = \int_{\Gamma_e} \psi \widetilde{(\varphi \mathbf{I}_d)}_N \cdot \hat{\mathbf{n}}\, d\Gamma_e - \int_{\Omega_e} \nabla \psi \cdot (\varphi \mathbf{I}_d)_N^{(e)}\, d\Omega_e\end{aligned} \tag{B.3}$$

where $\mathbf{I}_d$ is the rank-d identity matrix, where $d$ denotes the spatial dimension of the problem. With little loss in generality let us write the numerical flux or the first terms on the right-hand side of Eq. (B.3) in the following form

$$\begin{aligned}\widetilde{\mathbf{U}}_N = \{\{\mathbf{U_N}\}\} - \frac{\lambda}{2}[\![\varphi_N]\!] \\ \widetilde{(\varphi \mathbf{I}_d)}_N = \{\{\varphi_N \mathbf{I}_d\}\} - \frac{\lambda}{2}[\![\mathbf{U}_N]\!]\end{aligned} \tag{B.4}$$

where $\{\{\cdot\}\}$ and $[\![\cdot]\!]$ represent the average and jump operators, and $\lambda$ is the eigenvalue of the flux Jacobian. With this in mind, we can rewrite Eq. (B.3) in the following matrix form (see, e.g., [1])

$$\begin{aligned}\int_{\Omega_e} \psi \nabla \cdot \mathbf{U}_N \, d\Omega_e = \left(\boldsymbol{\mathcal{C}}^{(e,k)}\right)^{\mathcal{T}} \mathbf{U}_N - \mathcal{J}^{(e,k)} \varphi_N - \left(\widetilde{\boldsymbol{D}}^{(e)}\right)^{\mathcal{T}} \mathbf{U}_N \\ \int_{\Omega_e} \psi \nabla \varphi_N \, d\Omega_e = \boldsymbol{\mathcal{C}}^{(e,k)} \varphi_N - \mathcal{J}^{(e,k)} \mathbf{U}_N - \widetilde{\boldsymbol{D}}^{(e)} \varphi_N\end{aligned} \tag{B.5}$$

where $\boldsymbol{\mathcal{C}}^{(e,k)}$, $\mathcal{J}^{(e,k)}$, and $\widetilde{\boldsymbol{D}}^{(e)}$ represent the centered, jump, and weak form differentiation matrices. Note that $\boldsymbol{\mathcal{C}}^{(e,k)}$ and $\widetilde{\boldsymbol{D}}^{(e)}$ are vector matrices whereas $\mathcal{J}^{(e,k)}$ is a scalar matrix. Further note that the matrix superscript $(e, k)$ denotes that the matrices of the element $e$ also need the contribution of its neighbor $k$. Matrices with superscript $(e)$ are purely local element-wise quantities.

Taking all of this into account and assuming a backward Euler time-integration, we can now write Eq. (B.2) as follows

$$\begin{aligned}M^{(e)}\varphi + \Delta t \left[\left(\boldsymbol{\mathcal{C}}^{(e,k)}\right)^{\mathcal{T}} \mathbf{U} - \mathcal{J}^{(e,k)}\varphi - \left(\widetilde{\boldsymbol{D}}^{(e)}\right)^{\mathcal{T}} \mathbf{U}\right] = R_\varphi \\ M^{(e)}\mathbf{U} + \Delta t \left[\boldsymbol{\mathcal{C}}^{(e,k)}\varphi - \mathcal{J}^{(e,k)}\mathbf{U} - \widetilde{\boldsymbol{D}}^{(e)}\varphi\right] = R_{\mathbf{U}}\end{aligned} \tag{B.6}$$

where terms on the left-hand side are defined at the new time level $n + 1$ and terms on the right are at the current time level $n$; we also have dropped the subscript $N$ for convenience.

Letting

$$\begin{aligned}\widehat{M} = M^{(e)} - \Delta t \mathcal{J}^{(e,k)} \\ \widehat{\boldsymbol{D}} = \boldsymbol{\mathcal{C}}^{(e,k)} - \widetilde{\boldsymbol{D}}^{(e)}\end{aligned} \tag{B.7}$$

allows us to write Eq. (B.6) as follows

$$\begin{aligned}\widehat{M}\varphi + \Delta t \widehat{\boldsymbol{D}}^{\mathcal{T}} \mathbf{U} = R_\varphi \\ \widehat{M}\mathbf{U} + \Delta t \widehat{\boldsymbol{D}}\varphi = R_{\mathbf{U}}.\end{aligned} \tag{B.8}$$

Finally, applying a block LU decomposition yields the Schur complement

$$\widehat{M}\varphi - \Delta t^2 \widehat{\boldsymbol{D}}^{\mathcal{T}} \widehat{M}^{-1} \widehat{\boldsymbol{D}}\varphi = R_\varphi - \Delta t \widehat{\boldsymbol{D}} \widehat{M}^{-1} R_{\mathbf{U}}. \tag{B.9}$$

Looking at Eq. (B.7) and Eq. (B.9) we note that if $\mathcal{J}^{(e,k)}$ is not empty, then $\widehat{M}^{-1}$ is not block diagonal and thereby requires a global solution. However, if $\mathcal{J}^{(e,k)}$ is empty, which it will be for the continuous Galerkin method and for the DG method using centered fluxes, then $\widehat{M}^{-1}$ is block diagonal and trivial to invert - in fact, if we use inexact integration (as we do throughout this paper) it will be diagonal. Therefore, to make the Schur complement feasible requires using centered fluxes for the implicit linear operators. We say *feasible* because inverting $\widehat{M}$ will increase the cost of the Schur complement by a factor dependent on how often it appears in the Schur form. To give you a sense of what we mean, in Eq. (B.9), the inverse of $\widehat{M}$ appears once on the left and also on the right. Therefore, for non-block-diagonal $\widehat{M}$, Eq. (B.9) would require three global linear solves.

# C Discrete Formulations of the Schur Complement for the IMEX Discretizations

Defining the weak form differentiation matrix $\tilde{\boldsymbol{D}}^{(e)}$ as

$$\tilde{\boldsymbol{D}}^{(e)}_{ij} = \int_{\Omega_e} \nabla\psi_i \psi_j \, d\Omega_e,$$

the strong form differentiation matrix as

$$\boldsymbol{D}^{(e)}_{ij} = \int_{\Omega_e} \psi_i \nabla\psi_j \, d\Omega_e,$$

and the centered flux matrix $\boldsymbol{\mathcal{C}}^{(e,k)}$ as

$$\boldsymbol{\mathcal{C}}^{(e,k)}_{ij} = \int_{\Gamma_e} \psi_i \psi_j \hat{\mathbf{n}}^{(e,k)} \, d\Gamma_e,$$

for element $e$ with neighboring element $k$ with an outward pointing normal vector $\hat{\mathbf{n}}^{(e,k)}$ from $e$ to $k$. Then the weak variational form for $\nabla \cdot \mathbf{Q}$ can be written as

$$\int_{\Omega_e} \psi \nabla \cdot \mathbf{Q} \, d\Omega_e = \int_{\Gamma_e} \psi \mathbf{Q} \cdot \hat{\mathbf{n}} \, d\Gamma_e - \int_{\Omega_e} \nabla\psi \cdot \mathbf{Q} \, d\Omega_e = (\boldsymbol{\mathcal{C}}^{(e,k)})^{\mathcal{T}} \mathbf{Q} - (\tilde{\boldsymbol{D}}^{(e)})^{\mathcal{T}} \mathbf{Q} = \hat{\boldsymbol{D}}^{\mathcal{T}} \mathbf{Q}$$

where $\hat{\boldsymbol{D}} = \left(\boldsymbol{\mathcal{C}}^{(e,k)} - \tilde{\boldsymbol{D}}^{(e)}\right)$. Letting $\breve{\boldsymbol{D}} = \widehat{M}^{-1} \hat{\boldsymbol{D}}^{\mathcal{T}}$, where

$$\widehat{M}^e_{ij} = M^e_{ij} - \alpha \mathcal{J}^{(e,k)}_{ij}$$

is the augmented mass matrix and $\mathcal{J}^{(e,k)}$ is the jump term (that is zero for centered flux and continuous fields), $\widehat{M}^{-1} \nabla \cdot \mathbf{Q} = \breve{\boldsymbol{D}}^{\mathcal{T}} \mathbf{Q}$ and $\widehat{M}^{-1} \nabla \mathbf{Q} = \breve{\boldsymbol{D}} \mathbf{Q}$, the discrete weak Schur form for Set2C (Eq. (3.15)) can be obtained as

$$\begin{aligned} P_{tt} - \alpha^2 \frac{\gamma P_0}{\Theta_0} \breve{\boldsymbol{D}}^{\mathcal{T}} \left( \left[ h_{2C} \left\{ \mathbf{A}^{-1} \left( \breve{\boldsymbol{D}} P_{tt} + \frac{\breve{\boldsymbol{D}}\phi}{h_{2C}} \frac{\Theta_0 P_{tt}}{\gamma P_0} \right) \right\} \right] \right) \\ = \frac{\gamma P_0}{\Theta_0} \widehat{\Theta} - \alpha \frac{\gamma P_0}{\Theta_0} \breve{\boldsymbol{D}}^{\mathcal{T}} \left( \left[ h_{2C} \left\{ \mathbf{A}^{-1} \left( \widehat{\mathbf{U}} - \alpha \left[ \widehat{\rho} \, \breve{\boldsymbol{D}}\phi - \frac{\breve{\boldsymbol{D}}\phi}{h_{2C}} \widehat{\Theta} \right] \right) \right\} \right] \right). \end{aligned} \tag{C.1}$$

where $\alpha$ is the IMEX coefficient defined in Sec. 3.1. The discrete strong Schur form can be obtained by expanding the divergence operators in Eq. (3.15) using the product rule and substituting $\breve{\boldsymbol{D}}$ and $\breve{\boldsymbol{D}}^{\mathcal{T}}$ appropriately. Similarly, the discrete weak Schur form for Set3C (Eq. (3.26)) can be written as

$$
\begin{aligned}
&P_{tt} - \alpha^2 (\gamma - 1) \left( \breve{\boldsymbol{D}}^{\mathcal{T}} \left[ h_{3C} \left\{ \mathbf{A}^{-1} \left( \breve{\boldsymbol{D}} P_{tt} + \frac{\breve{\boldsymbol{D}} \phi}{(\gamma - 1)(h_{3C} - \phi)} P_{tt} \right) \right\} \right] \right) \\
&\quad - \alpha^2 (\gamma - 1) \phi \breve{\boldsymbol{D}}^{\mathcal{T}} \left\{ \mathbf{A}^{-1} \left( \breve{\boldsymbol{D}} P_{tt} + \frac{\breve{\boldsymbol{D}} \phi}{(\gamma - 1)(h_{3C} - \phi)} P_{tt} \right) \right\} \\
&\quad = (\gamma - 1) \left( \widehat{E} - \phi \widehat{\rho} \right) \\
&\quad - \alpha (\gamma - 1) \left( \breve{\boldsymbol{D}}^{\mathcal{T}} \left[ h_{3C} \left\{ \mathbf{A}^{-1} \left( \widehat{\mathbf{U}} - \alpha \frac{\breve{\boldsymbol{D}} \phi}{(h_{3C} - \phi)} \left( h_{3C} \widehat{\rho} - \widehat{E} \right) \right) \right\} \right] \right) \\
&\quad - \alpha (\gamma - 1) \phi \breve{\boldsymbol{D}}^{\mathcal{T}} \left\{ \mathbf{A}^{-1} \left( \widehat{\mathbf{U}} - \alpha \frac{\breve{\boldsymbol{D}} \phi}{(h_{3C} - \phi)} \left( h_{3C} \widehat{\rho} - \widehat{E} \right) \right) \right\}
\end{aligned}
\tag{C.2}
$$

and the discrete strong Schur form can be obtained by substituting $\breve{\boldsymbol{D}}$ and $\breve{\boldsymbol{D}}^{\mathcal{T}}$ appropriately in Eq. (3.25). Note that the strong forms require two applications of integration by parts as described in [15, 53].